\documentclass[12pt]{article}
\usepackage{amsmath,amssymb,bbold,color,xcolor,graphicx,float,mathrsfs, braket, array, makecell}
\usepackage{tikz}
\usetikzlibrary{knots, intersections}
\usepackage{ytableau}
\usepackage[colorlinks=true,linkcolor=blue]{hyperref}
\newcommand{\RNum}[1]{\uppercase\expandafter{\romannumeral #1\relax}}
\usepackage[autostyle]{csquotes}
\usepackage[titles]{tocloft}

\hypersetup{pdfborder=0 0 0}
\begin{document}
\begin{titlepage}
\begin{center}
\vspace{10mm}
{\LARGE \textbf{A survey of knots and quivers}}
\vspace{12mm}

\renewcommand\thefootnote{\mbox{$\fnsymbol{footnote}$}}
\large{Shivrat Sachdeva}\footnote{shivrats@imsc.res.in}

\vspace{10mm}

\end{center}
\vspace{12mm}

\noindent
This survey explores knot polynomials and their categorification, culminating in the homological invariants of knots. We begin with an overview of classical knot polynomials, progressing towards the superpolynomial and its role in unifying various knot homologies. Along the way, we provide physical and geometric insights into the unification of the \(sl(N)\) Khovanov-Rozansky and the knot Floer homology. We then turn our attention to the intriguing correspondence between knots and quivers, examining how this perspective sheds light on the integrality of BPS states encoded in the Labastida-Mariño-Ooguri-Vafa (LMOV) invariants. We will further investigate the knot-quiver correspondence from a physics and geometric side and study the 3d \(\mathcal{N}=2\) theory \(T[Q_K]\) for the quivers. 

\end{titlepage}
\setcounter{footnote}{0}
\renewcommand\thefootnote{\mbox{\arabic{footnote}}}

\hrule
\tableofcontents
\bigskip
\hrule

\addtolength{\parskip}{8pt}
\section{Introduction}
In recent times there has been a remarkable surge towards the union between quantum field theory, string theory, topology and knot theory, sparked by the seminal paper by Witten \cite{Witten:1988hf}. Following this breakthrough, other important mathematical developments followed such as polynomial knot invariants and their categorification, volume conjectures, \(A-\)polynomials and much more, which further have been interpreted through the lens of physics. 
\\ \\
In Section \ref{KnotInvariants}, we present a brief introduction to polynomial knot invariants. We then introduce the concept of categorification, in the simplest words it means: Given some data, we can encode it in a more sophisticated object (categorified data) and ask how to extract the original data from this more complicated data set. The most common example of this is the Khovanov homology \cite{khovanov1999categorificationjonespolynomial}, which produces the Jones polynomial as its Euler characteristic. The homology can be viewed as this sophisticated data set, from which the Euler characteristic can be extracted. Categorification is much more powerful than its decategorified counterpart, thus is useful in classifying manifolds and topological spaces. 
\\ \\
In Section \ref{Knotsuperpolynomial} we review the framework presented in the work \cite{dunfield2005superpolynomialknothomologies}, unifying the \(sl(N)\) Khovanov-Rozansky homology with the knot Floer homology using a triply-graded homology. We start off by introducing the Poincaré (super)polynomial for the homology \(\mathscr{H}_{i,j,k}\). Furthermore, we explain in great detail the idea of color graded differentials, which are a part of the package of the tri-graded homology. We will review the techniques and explicitly work out several examples, taking the homology with respect to the colored differentials of the chain complex. In the end, we shed light on the unification procedure with the knot Floer homology and paint a physical picture of the triply-graded homology and its role in the study of Gromov-Witten (GW) theory and BPS spectrum. 
\\ \\
In Section \ref{Knots&Quivers} we present a connection between knots, generated by brane engineering on the string theory side and the representation theory of quivers \cite{Kucharski_2017,Kucharski_2019}. Topological string amplitudes encode the Gromov-Witten invariants, which from the physics perspective capture the BPS states. On the closed string theory side, these invariants are also refereed to as the Gopakumar-Vafa invariants (closed BPS invariants) \cite{gopakumar1998mtheorytopologicalstringsii, Gopakumar:1998ii} and are further related to the Donaldson-Thomas invariants for the ambient Calabi-Yau manifold. With additional branes wrapping the Lagrangian conormal, the open string amplitudes capture the colored HOMFLY-PT polynomial \cite{Ooguri_2000}. The corresponding analog on the open string side are the BPS invariants also called the Labastida-Mariño-Ooguri-Vafa (LMOV) invariants \cite{Ooguri_2000, Labastida_2001, Labastida_2000}. We show that using the knot-quiver correspondence, the integrality of the LMOV invariants for symmetric quivers can be proven. We highlight the fact that the map between knots and quivers is not an bijective one. Finally, we illustrate the role of choice of knot normalization factor, and its consequences on the quiver side. 
\\ \\
In Section \ref{Geometry_KQ} we present the results and claims made in \cite{ekholm2020physicsgeometryknotsquiverscorrespondence}. We study the knot-quiver correspondence from a two-fold perspective; on one hand we have topological string theory and on the other lies holomorphic curves and sympletic geometry. We begin by presenting the M-theory construction of knots and quivers, and then give present the dual picture in terms of holomorphic disks ending on the Lagrangian conormal \(\mathscr{L}_K \subset X\). Then we move onto introducing the quiver gauge theory \(T[Q_K]\) associated to a quiver \(Q_K\) and show the relationship of the quiver superpotential to the Gromov-Witten disk potential. Furthermore, we will present a connection between the quiver data to the (self)linking and 4-chain intersection. Finally, we work out the example of an unknot and illustrate the results of the authors \cite{ekholm2020physicsgeometryknotsquiverscorrespondence}.
\\ \\
\bigskip
\goodbreak
\centerline{\bf Acknowledgments}
\noindent
I am grateful to Piotr Kucharski for his valuable inputs towards the draft and for the lively discussions that taught me the ideas and results presented in this work.
\section{Knot Invariants} \label{KnotInvariants}
Over the past few decades, we have seen a surge in the study of intersection between knots, strings and low-dimensional topology. Ever since the discovery of the Jones polynomial in the context of Chern-Simons theory \cite{Witten:1988hf}, extensive research has been conducted in studying quantum invariants of knots and links. 
\\ \\
We begin by defining an oriented knot colored by the fundamental representation of the Lie group \(sl(N)\), which leads to the quantum \(sl(N)\) invariant \(P^K_N(q)\), and defines the HOMFLY-PT polynomial for us. For a knot \(K\) embedded in \(S^3\), its normalized HOMFLY-PT polynomial \(P^K(a=q^N,q)\) is described using the skein relation,
\begin{equation} 
    q^N P^K_{N}\left(
    \begin{tikzpicture}[baseline={([yshift=-.5ex]current bounding box.center)}, scale=0.7]
        \begin{knot}[clip width=3pt, flip crossing=0]
            \strand[thick, ->] (-0.4,-0.4) to (0.4,0.4);
            \strand[thick, ->] (0.4,-0.4) to (-0.4,0.4);
        \end{knot}
    \end{tikzpicture}
    \right)-q^{-N} P^K_{N}\left(
    \begin{tikzpicture}[baseline={([yshift=-.5ex]current bounding box.center)}, scale=0.7]
        \begin{knot}[clip width=3pt, flip crossing=1]
            \strand[thick, ->] (-0.4,-0.4) to (0.4,0.4);
            \strand[thick, ->] (0.4,-0.4) to (-0.4,0.4);
        \end{knot}
    \end{tikzpicture}
    \right)  = 
    (q- q^{-1}) P^K_{N}\left(
    \begin{tikzpicture}[baseline={([yshift=-.5ex]current bounding box.center)}, scale=0.7]
        \draw[thick, ->] (-0.4,-0.4) to[out=70, in=290] (-0.3,0.4);
        \draw[thick, ->] (0.4,-0.4) to[out=110, in=250] (0.3,0.4);
    \end{tikzpicture}
    \right)
    \end{equation}
To switch between the normalized to the unnormalized polynomial, we divide the unnormalized knot polynomial by the unnormalized unknot factor given by \(\overline{P}\left(
    \begin{tikzpicture}[baseline={([yshift=-.5ex]current bounding box.center)}, scale=0.7]
        \draw[thick, ->] (0,0) circle (0.3);
        \draw[thick, ->] (0.3,0) arc (0:180:0.3);
    \end{tikzpicture}
    \right) =(a-a^{-1})/(q-q^{-1})\). The HOMFLY-PT polynomial describes the \(sl(N)\) knot invariant for all \(N\) and the Jones polynomial \(J_K\) is simply \(P_2^K(a=q^2,q)\). Furthermore, taking the limit \(a=1\), in the normalized HOMFLY-PT polynomial, gives the Alexander polynomial \(\Delta^K(a=1,q)\). 
\\ \\
A major drawback of these such knot invariants is that they cannot distinguish between mutants \cite{doi:10.1142/S0218216596000163}, however their respective categorifications can. The general idea for a knot \(K\) is to construct a bi-graded homology \(H_{i,j}(K)\) whose graded Euler characteristic with respected to a grading gives a knot polynomial.  
\subsection{Categorification}
Categorification for a 3-dimensional TQFT manifests a 4-dimensional TQFT, from which a 3D theory can be extracted by dimensional reduction \cite{Crane_1994,Gukov:2007ck}. Each TQFT associates to itself a class of geometric objects,
\begin{table}[H]
    \centering
    \renewcommand{\arraystretch}{1.5}
    \begin{tabular}{ |c|c|c| } 
        \hline
        \textbf{Geometry} & \textbf{3D TQFT} & \textbf{4D TQFT} \\
        \hline
        \makecell{3-manifold \(M\),\\ embedded knot \(K \subset M\)} 
        & Poincaré polynomial \(P^K\) 
        & vector space \(\mathscr{H}_K\) \\
        \hline
        2-manifold \(\Sigma\) 
        & vector space \(\mathscr{H}_K\) 
        & category \(\mathrm{Cat}_{\Sigma}\) \\
        \hline
    \end{tabular}
    \caption{Embedded structures in 3D and 4D TQFT.}
    \label{tab:tqft_table}
\end{table}
To any geometric object of a given dimension, a categorified TQFT assigns an object with one higher dimension than its decategorified dual. Some famous examples are the Donaldson theory \cite{Witten:1988ze} and the Seiberg-Witten theory \cite{witten1994monopolesfourmanifolds}. 
\\ \\
Khovanov in his paper \cite{khovanov1999categorificationjonespolynomial} constructed the categorified version of the Jones polynomial. In his work, he associates a chain complex to a knot \(K\), with the homology of this chain complex being invariant under the Reidemeister moves, as a consequence it is also an invariant for the knot \(K\). The Khovanov homology \(H^{Kh}_{i,j}(K)\) is a bi-graded homology whose Euler characteristic gives back the Jones polynomial,
\begin{equation}
  J^K(q)=\sum_{i,j}(-1)^iq^j \text{dim}\ H^{Kh}_{i,j}(K)  
\end{equation}
where \(i-\) is the Jones grading and \(j-\) is the homological grading. Khovanov's theory was later generalized to categorify all \(sl(N)\) quantum invariants. The Khovanov-Rozansky homology \cite{khovanov2004matrixfactorizationslinkhomology} is a doubly-graded homology \(HRK_{i,j}^{N}(K)\) whose graded Euler characteristic gives the \(sl(N)\) polynomial invariant \(P^K_N\),
\begin{equation}
  P^K_N(q)=\sum_{i,j} (-1)^j q^i \text{dim}HKR_{i,j}^{N}(K) 
\end{equation} \label{HKR}
For \(N=2\), we retrieve the original Khovanov homology. An equivalent categorified homology theory also exists for the Alexander polynomial, given by
\begin{equation}
    \Delta^K(q)=\sum_{i,j} (-1)^j q^i \text{dim} \widehat{HFK}_j(K;i)
\end{equation}
where \(\widehat{HFK}_j(K,i)\) is the knot Floer homology. One should keep in mind that knot homologies as well as the colored differentials \(\{d_N\}\) admit a physical description in terms of refined BPS invariants \cite{Gukov_2005, Gukov_2016, Gukov_2013}. 
\section{Knot superpolynomial} \label{Knotsuperpolynomial}
The knot superpolynomial was introduced as a mean to unify the Khovanov-Rozansky \(sl(N)\) homology and the knot Floer homology \cite{dunfield2005superpolynomialknothomologies}. The proposal goes as follows: There must exist a finite superpolynomial of the form \(\mathscr{P}^K \in \mathbb{Z}_{\geq0}[a^{\pm1},q^{\pm1},t^{\pm1}]\) such that,
\begin{equation}
    P^K(a,q)=\mathscr{P}^K(a,q,t=-1)
\end{equation}
and for the large \(N\) limit, we get the following expression,
\begin{equation}
    KhR_N(q,t)=\mathscr{P}^K(a=q^N,q,t)
\end{equation}
One views this superpolynomial \(\mathscr{P}^K\) as the Poincaré polynomial of a triply-graded homology \(\mathscr{H}_{i,j,k}(K)\) which categorifies the HOMFLY-PT polynomial. The idea is that one can obtain the bi-graded \(sl(N)\) homology by taking the homology of triply-graded complex \(\mathscr{H}_*(K)\) with respected to differential \(\{d_N \}\). 
\begin{figure}[H]
    \centering
    \includegraphics[width=0.95\linewidth]{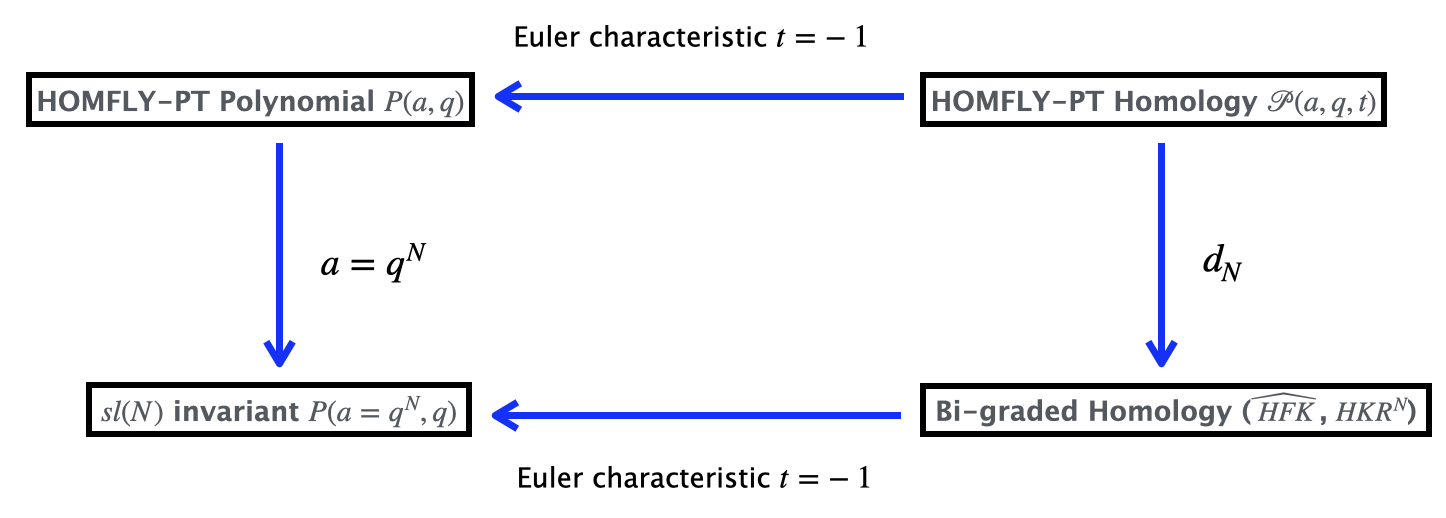}
    \caption{This figure depicts the various gradings and their connections.}
    \label{fig:Hom diagram}
\end{figure}
The Poincaré polynomial for the triply-graded homology is given by,
\begin{equation}
    \mathscr{P}^K(a,q,t)=\sum_{i,j,k}a^iq^jt^k\text{dim}\mathscr{H}_{i,j,k}(K)
\end{equation}
where the homology \(\mathscr{H}_*(K)\) equipped with a family colored differentials. For \(N > 0\), the homology \((\mathscr{H}_*(K),d_N)\) is isomorphic to the \(sl(N)\) Khovanov-Rozansky homology, whereas for \(N=0\), \((\mathscr{H}_*(K),d_0)\) is isomorphic to the knot Floer homology. 
\\ \\
The colored differentials are defined in the following way,
\begin{enumerate}
    \item \(N>0\), \(d_N\) is triply-graded under \((-2,2N,-1)\)
    \begin{equation}
        d_N : \mathscr{H}_{i,j,k}(K) \rightarrow \mathscr{H}_{i-2,j+2N,k-1}(K)
    \end{equation}
    \item \(N<0\), \(d_N\) is triply-graded under \((-2,2N,-1+2N)\)
     \item \(N=0\), \(d_0\) is triply-graded under \((-2,0,-3)\)
    \item \(d_N\)'s are anticommuting, i.e. \(d_N^2=0\) for \(N \in \mathbb{Z}\)
    \item There is an involution symmetry \(\phi\) s.t 
    \begin{align}
        & \phi: \mathscr{H}_{i,j,*}(K) \rightarrow \mathscr{H}_{i,-j,*}(K) \nonumber \\ &\phi d_N=d_{-N}\phi \;\ \text{for} \;\ N \in \mathbb{Z}
    \end{align}
\end{enumerate}
Let's illustrate the action of the colored differentials on corresponding superpolynomial of the some simple knots.
\subsection{Trefoil Knot \texorpdfstring{$T_{2,3}$}{T(2,3)}}
Let's illustrate how to compute the homology with respected to the colored differentials for the trefoil knot and later on move onto a more non-trivial example of the torus knot \(T_{3,4}\). 
\\ \\
\textit{Trefoil knot \(3_1\):} The superpolynomial for the negative trefoil knot \(T_{2,3}\).
\begin{equation}
    \mathscr{P}(T_{2,3})=a^2q^2t^0+a^2q^2t^2+a^4a^0t^3 \label{3_1 sp}
\end{equation}
Expressed in terms of the homological diagram,
\begin{figure}[H]
    \centering
    \includegraphics[width=0.30\linewidth]{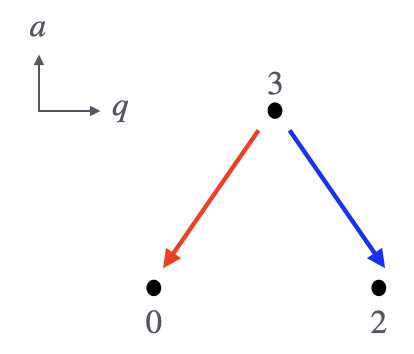}
    \caption{Homological diagram for the trefoil with the non-trivial colored differentials \(d_{\pm1}\).}
    \label{fig:3_1 trefoil}
\end{figure}
For the trefoil, there are only two non-zero differentials \(d_{1}\) and \(d_{-1}\). For \(N>1\), all the differentials act trivially on the superpolynomial, giving us \(sl(N)\) invariants simply by substituting \(a=q^{N}\) in \eqref{3_1 sp} for all \(N >1\). For \(N=1\), \(d_1\) kills the right-hand generator, giving \(\mathscr{P}_1(T_{2,3})=1=q^0t^0\), which is simply the \(sl(1)\) invariant. One should note that the  \(sl(1)\) invariant is always trivial for any knot (monomial). 
\subsection{Torus Knot  \texorpdfstring{$T_{3,4}$}{T(3,4)}}
Here we do a similar homology computing exercise, but for a more complex knot. \\ \\
\textit{Torus knot \(T_{3,4}\):} The superpolynomial for the knot \(T_{3,4}\) is given by,
\begin{equation}
     \mathscr{P}(T_{3,4})=a^{10}t^8+a^8(q^{-4}t^3+q^{-2}t^5+t^5+q^2t^7+q^4t^7)+a^6(q^{-6}+q^{-2}t^2+t^4+q^2t^4+q^6t^6)
\end{equation}
Here we have five non-zero colored differentials \(d_{1},d_{2},d_{-1},d_{-2}\) and \(d_0\). 
\begin{figure}[H]
    \centering
    \includegraphics[width=0.60\linewidth]{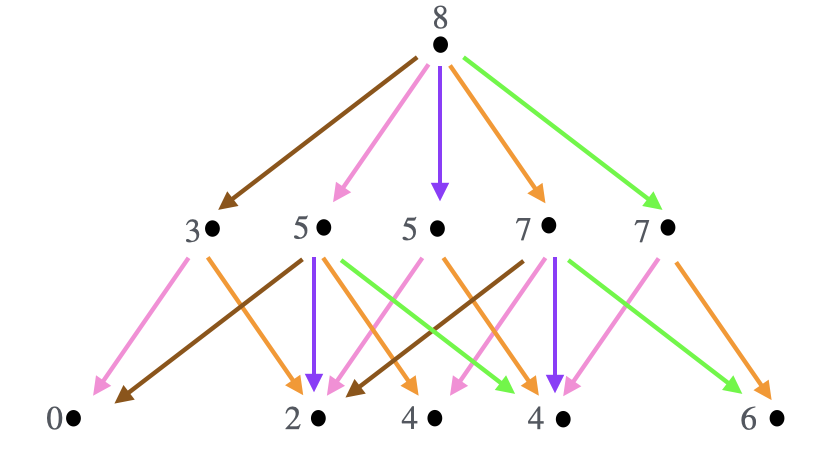}
    \caption{Homological diagram for the knot \(T_{3,4}\). The bottommost row has \(a-\)grading 6, the leftmost dot has \(q-\)grading -6 and the vertical axis of symmetry lies at \(q-\)grading 0. }
    \label{fig:T34 knot}
\end{figure}
By taking the differential with respected to \(d_2\) of \(\mathscr{H}_*\) and substituting \(a=q^2\), we get the Khovanov homology ,
\begin{equation}
\mathscr{P}_2(T_{3,4})=q^6+q^{10}t^2+q^{12}t^3+q^{12}t^4+q^{16}t^5
\end{equation}
and as before \(\mathscr{P}_1(T_{3,4})=1\), the bottom leftmost dot survives.
\subsection{Connection to knot Floer homology}
In order to get the knot Floer homology, we introduce a new homological grading on \(\mathscr{H}(K)\) defined by \(t'(x)=t(x)-a(x)\). The knot Floer homology is defined by the action of the \(N=0\) differential with respect to the new homological grading \(t'\). The Poincaré polynomial is defined by,
\begin{equation}
    HFK(q,t)=\sum_{i,j}q^it^j\text{dim}\widehat{HFK}_j(K;i) \label{hfk}
\end{equation}
After the new grading, the superpolynomial becomes,
\begin{equation}
    \mathscr{P'}^K(a,q,t)=\mathscr{P}^K(a=at^{-1},q,t)
\end{equation}
The differential \(d_0\) reduces the \(t'-\)grading by 1. Substituting \(a=1\) and taking the grading with respected to \(d_0\), we retrieve \(\mathscr{P}_0^K(q,t)\), which is Poincaré polynomial of the homology which categorifies the Alexander polynomial \(\Delta^K(q)=P^K(a=1,q)\). Moreover, \(\mathscr{P}_0^K=HFK(q,t)\), where \(HFK\) is the Poincaré polynomial for the knot Floer homology as defined in \eqref{hfk}.
\begin{figure}[H]
    \centering
    \includegraphics[width=0.30\linewidth]{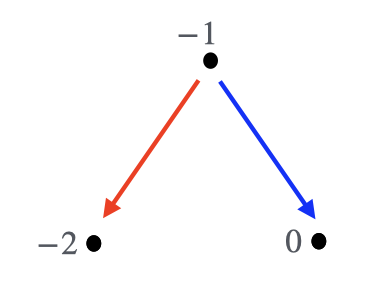}
    \caption{Homological diagram for the trefoil knot with respect to the new homological grading \(t'\).}
    \label{fig:floer}
\end{figure}
In the case of the trefoil knot, depicted in Figure \ref{fig:floer}, the differential \(d_0\) acts trivially, 
\begin{equation}
    \mathscr{P}_0(T_{2,3})=\mathscr{P}(T_{2,3})(a=t^{-1},q,t)=q^{-2}t^{-2}+t^{-1}+q^2
\end{equation}
which is the Poincaré polynomial \(HFK(T_{2,3})\) of the knot Floer homology. Doing a same exercise for the knot \(T_{3,4}\) yields
\begin{equation}
    \mathscr{P}_0(T_{3,4})=q^{-6}t^{-6}+q^{-4}t^{-5}+t^{-2}+q^4t^{-1}+q^6
\end{equation}
agrees with \(HFK(T_{3,4})\).
\subsection{A Geometric picture}
It turns out that the triply-graded homology of the unreduced knot is closely related to the geometry of holomorphic curves. The physical setup required to understand this relation is a non-compact Calabi-Yau 3-fold \(X\) and a Lagrangian submanifold \(\mathscr{L}_K \subset X\) \cite{Taubes_2006}, i.e. the knot conormal \(\mathscr{L}_K\). For every embedded knot \(K \subset S^3\), we define a \(\mathscr{L}_K \; \text{and} \; X\). The Calabi-Yau manifold is chosen to be the resolved conifold \cite{Ooguri_2000},
\begin{equation}
    \mathscr{O}(-1) \oplus \mathscr{O}(-1) \rightarrow \mathbb{CP}^1
\end{equation}
We can study the embedded holomorphic Riemann surfaces in \(X\) with boundaries ending on the knot conormal \(\mathscr{L}_K\), 
\begin{equation}
    (\Sigma, \partial\Sigma) \rightarrow (X,\mathscr{L}_K)
\end{equation}
Given such a class of embedded surfaces satisfying certain conditions as stated in \cite{dunfield2005superpolynomialknothomologies}, we define \(\mathscr{M}_{g,Q}(X, \mathscr{L}_K)\) to be the moduli space of embedded Riemann surfaces \(\Sigma\) with flat connections modulo gauge transformation. The cohomology groups of this moduli space \(H^k(\mathscr{M}_{g,Q})\) is described using three integers: \(k\) is the degree, \(g\) the genus and \(Q \in H_2(X, \mathscr{L}_K; \mathbb{Z}) \cong \mathbb{Z}\) is the relative homology class; which also turn out to be the gradings of our triply-graded homology. 
\\ \\
The Euler characteristic \(\chi(\mathscr{M}_{g,Q})\) of the cohomology group \(H^k(\mathscr{M}_{g,Q})\) encodes the information about the integer BPS invariants, keeping count of all-genus open Gromov-Witten invariants \cite{Ooguri_2000, Labastida_2000}. This turns out to be equivalent to the Euler characteristic of the triply-graded homology \(\mathscr{H}_*(K)\). This equivalence enables us to encode all the information about the GW invariants into a finite-set of non-zero BPS integers.  
\section{Knots and Quivers} \label{Knots&Quivers}
The knots-quiver correspondence provides a relation between two distinct fields of mathematics, knot theory and the representation theory of quivers. Precisely speaking, various knot invariants are encoded in the moduli space of quivers representations. Let's try to understand both sides of the coins first, before connecting them via the knot-quiver correspondence.
\subsection{Knot and TQFT}
In 1988, Witten in a revolutionary paper \cite{Witten:1988hf} discovered that the Jones polynomial originates from the Chern-Simons (CS) gauge theory in 3-dimensions. The Chern-Simons theory is topological by nature and the classical action on a manifold \(M\) with an associated gauge group \(G\) with a parameter \(k\) (inverse of the coupling constant) is given by,
\begin{equation}
    S_{CS}=\frac{k}{4\pi}\int_M\text{Tr}\left(A \wedge A +\frac{2}{3}A \wedge A \wedge A \right)
\end{equation}
where the trace \(\text{Tr}\) is over the representations of the gauge group \(G\) and the parameter \(k\) is required to an integer, so that the integrand is invariant under gauge transformations, whose non-trivial elements correspond to \(\pi_3(G) =\mathbb{Z}\). 
\\ \\
Witten used the topological invariance of the CS-theory on an arbitrary 3-manifold \(M\) with embedded knots, by cutting and gluing the 3-manifold into pieces, solving the theory piecewise and then gluing the pieces back together. A precise connection between knot and gauge theory and beyond can be found in Mina Aganagic's article \cite{Aganagic:2015tka}.
\begin{figure}
    \centering
    \includegraphics[width=0.5\linewidth]{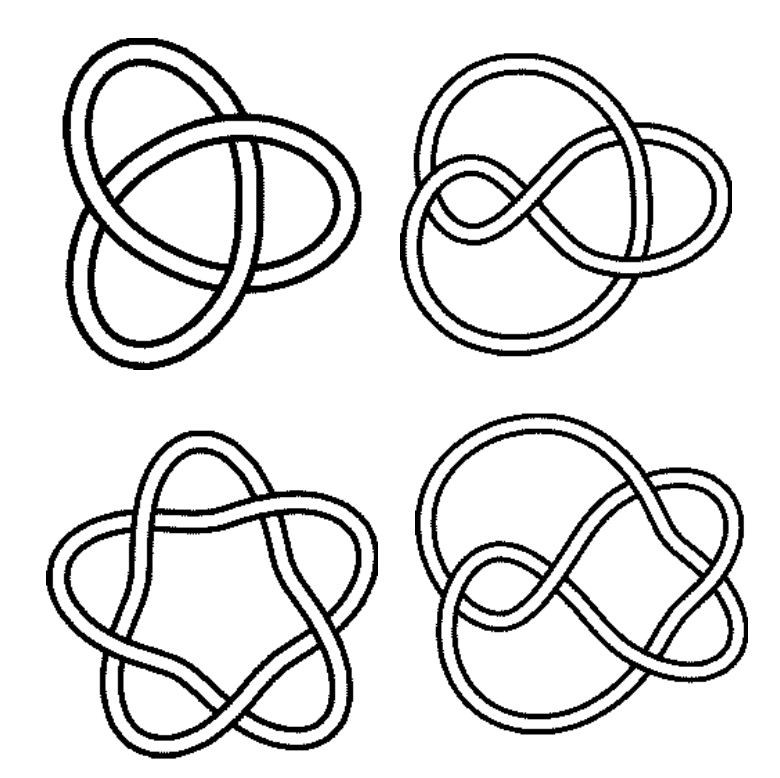}
    \caption{Knots \(3_1\), \(4_1\), \(5_1\) and \(5_2\). Images taken from the Knot Atlas.}
    \label{fig:knot_atlas}
\end{figure}
We are interested in writing the knot HOMFLY-PT polynomial using the CS-theory. Suppose there is a \(U(N)\) gauge connection \(A\) and the ambient three-manifold is \(M=S^3\). Then the unreduced HOMFLY-PT polynomial colored by a representation \(r\) of the gauge group \(U(N)\) is defined using the holonomy of the gauge connection \(A\) along the knot \(K \subset S^3 \),
\begin{equation}
    \overline{P}^K_r(a,q)=\int \mathcal{D}A\ e^{iS_{CS}[A]}\text{Tr}_r\ \text{Hol}_{\gamma}(A)
\end{equation}
where \(\gamma\) is any closed loop in \(S^3\) and the parameters \(a\) and \(q\) are defined as,
\begin{equation}
    q=\text{exp}\left(\frac{i\pi}{k+N}\right) \;\ a=\text{exp}\left(\frac{i\pi N}{k+N}\right)
\end{equation}
Restricting to symmetric representations \(S^r\), having Young diagram,
\[S^r=
\underbrace{
\ytableausetup{centertableaux}
\begin{ytableau}
~ & ~ & ~ & \none[\dots] & ~
\end{ytableau}
}_{\text{$r$ boxes}}
\]
Knot polynomials further exhibit an interesting property that will later become crucial in understanding the knot-quiver correspondence. Certain class of knots exhibit “exponential growth property", implying that the Poincaré superpolynomial of the \(r-\) colored (reduced) HOMFLY-PT homology obeys,
\begin{equation}
    \mathscr{P}^K_r(a,q,t) \simeq \left(\mathscr{P}^K_{\square}(a,q,t)\right)^n \label{expo_sp}
\end{equation}
where \(\square\) represents the fundamental representation \(S^1\) and \(\simeq\) denotes a bijection mapping between the monic terms of both sides. Note, that this property also holds for the unreduced homology generators. However in that case, one needs to take into consideration the exponential growth of the unknot factor as well. The relation \eqref{expo_sp}, becomes an equality in the semi-classical limit \(\hbar \to 0\),
\begin{equation}
    \mathscr{P}^K_r(a,1,t) = \left(\mathscr{P}^K_{\square}(a,1,t)\right)^n 
\end{equation}
The simplest knot that does not satisfy this property is \(9_{42}\). The proof, and its consequences in the context of the knot-quiver correspondence is discussed in \cite{ekholm2020physicsgeometryknotsquiverscorrespondence}.
\subsection{Quiver picture}
According to \cite{Ooguri_2000}, knots \(K \subset S^3\) are realized in topological string theory by wrapping \(N'\) A-branes on the on the Lagrangian knot conormal \(\mathscr{L}_K \subset T^*S^3\), shifted off the zero-section of \(S^3\) with additional \(N\) A-branes wrapping the base \(S^3\) in the deformed conifold \(T^* S^3\). After the large \(N\) transition, the deformed conifold transitions into the resolved conifold, base \(S^3\) shrinks to a point and we are only left with the branes wrapped around the \(\mathscr{L}_K\), encoding the information about knots. In the M-theory picture, it was argued that the colored HOMFLY-PT polynomials are encoded in the LMOV invariants, which also count the BPS spectrum of M2-branes ending on the M5-branes wrapping \(\mathscr{L}_K\). We will come back to this geometric point of view of knots and quivers at a later stage in Section \ref{Geometry_KQ}. 
\begin{table}[H]
    \centering
    \begin{tabular}{|c|c|}
    \hline
     \textbf{Knots}& \textbf{Quivers}\\
     \hline
      Homological Degrees, framing & Number of loops \\
      Colored HOMFLY-PT generators & Motivic generating series \\
      LMOV invariants & Motivic DT invariants \\
      Classical LMOV invariants & Numerical DT invariants \\
      Algebra of BPS states & Cohomological Hall Algebra \\
      \hline
    \end{tabular}
    \caption{Matching of quantities under the knot-quiver correspondence.}
    \label{KQ_table}
\end{table}
We start off by introducing the generating series of the \(S^r\) colored HOMFLY-PT polynomial using the Ooguri-Vafa generating function \cite{Ooguri_2000, labastida2001newpointviewtheory, Labastida_2001, Labastida_2000}
\begin{equation}
    Z(U,V)=\sum_r \text{Tr}_r U \ \text{Tr}_r V=\text{exp}\left(\sum_{n=1}^{\infty} \frac{1}{n}\text{Tr}{U^n}\text{Tr}{V^n} \right) \label{ov}
\end{equation}
where \(U =P \;\ \text{exp} \oint_K A \) is the holonomy of the \(U(N)\) Chern-Simons gauge field around the knot \(K\), and \(V\) is the source with the sum running over all representations \(r\). The LMOV conjecture states that the expectation value of \eqref{ov} is
\begin{align}
    \braket{Z(U,V)} &=\sum_r \overline{P}(a,q) \text{Tr}_r V \nonumber \\ &=\text{exp}\left(\sum_{n=1}^{\infty}\sum_r \frac{1}{n} f_r(a^n,q^n)\text{Tr}_r{V^n} \right)
\end{align}
where the expectation value of the \(U(N)\) holonomy is the unreduced HOMFLY-PT for a knot \(K\) with \(\braket{\text{Tr}_r U}=\overline{P}_r(a,q)=P_r^{0_1}P_r(a,q)\), where \(P_r^{0_1}=\frac{(a^2;q^2)_r}{(q^2;q^2)_r}\) is the unknot factor and the function \(f\) is identified as,
\begin{equation}
    f_r(a,q)=\sum_{i,j}\frac{N_{r,i,j}a^iq^j}{q-q^{-1}}
\end{equation}
encodes the LMOV invariants \(N_{r,i,j}\). Consider a one-dimensional source \(V=x\), thus \(\text{Tr}_r V \neq 0\) only for symmetric representations \(r=S^r\) and gives \(\text{Tr}_{S^r}\ x=x^r\). Putting these values back in the Ooguri-Vafa series \eqref{ov}, we get
\begin{align}
     P(x)=\sum_{r=0}^{\infty}\overline{P}_r(a,q)x^r=e^{\sum_{r, n\geq 1}\frac{1} {n}f_r(a^n,q^n)x^{nr}}  \label{p(x)}
\end{align}
where \(P(x)=\braket{Z(U,x)}\) and \(\overline{P}_r(a,q)\) denotes the \(S^r\) colored unreduced HOMFLY-PT polynomial for a knot \(K\). The LMOV function for the symmetric representation takes the form, 
\begin{equation}
    f_r(a,q)=f_{S^r}(a,q)=\sum_{i,j}\frac{N_{S^r,i,j}a^iq^j}{q-q^{-1}} 
\end{equation}
Now the LMOV invariants \(N_{S^r,i,j}\) are polynomials in the unreduced HOMFLY-PT polynomial \(\overline{P}_{d_1}(a^{d_2},q^{d_2})\). Using the properties of the Plethystic exponential, we can rewrite \eqref{p(x)} as,
\begin{equation}
    P(x) = \prod_{r \geq1,i,j,k \geq 0}(1-x^ra^iq^{j+2k+1})^{N_{S^r,i,j}} \label{pleth}
\end{equation}
In the limit \(\hbar \to 0\), we get the “classical" LMOV invariants, defined using the following ratio,
\begin{equation}
    y(x,a)=\lim_{q \to 1} \frac{P(qx)}{P(x)}=\lim_{q \to 1}\prod_{r \geq1,i,j,k \geq 0}\left(\frac{1-x^ra^iq^{r+j+2k+1}}{1-x^ra^iq^{j+2k+1}}\right)^{N_{S^r,i,j}}=\prod_{r \geq1;i \in \mathbb{Z}}(1-x^ra^i)^{-rb_{r,i}/2}
\end{equation}
where the “classical" LMOV invariants are defined to be,
\begin{equation}
    b_{r,i}=\sum_j N_{r,i,j} \label{classlmov}
\end{equation}
It turns out that the algebraic curve \(y=y(x,a)\) satisfies the zero locus of the \(A-\)polynomial \cite{Garoufalidis_2016, Cooper1994, aganagic2012largendualitymirror},
\begin{equation}
    A(x,y)=0
\end{equation}
A more general two-parameter deformed version of the \(A-\)polynomial was studied in \cite{Fuji:2013rra, Fuji_2013}. The claim is that this “super" \(A-\)polynomial encodes the color dependence of the HOMFLY-PT homology and the asymptotics of \(S^r\) colored homology in the large color limit \(r \to \infty\), followed by semi-classical limit \(q \to 1\).
\\ \\
Coming back to the KQ-correspondence, the statement is that the LMOV invariants can be written as a linear combinations of motivic Donaldson-Thomas (DT) invariants \(\Omega^{Q_K}_{d,s}\),
\begin{equation}
    \Omega^{Q_K}_{d,s}=\sum_{\substack{d_1,..,d_m \geq 0 \\ s \in \mathbb{Z}}}\Omega_{d_1...,d_m;s}q^sx_1^{d_1}...x_m^{d_m}
\end{equation}
where \(\Omega_{d_1...,d_m}\) are the numerical Donaldson-Thomas invariants which have been proven to be non-negative integers \(\mathbb{Z_{\geq0}}\) for any symmetric quiver \cite{Efimov_2012, kontsevich2008stabilitystructuresmotivicdonaldsonthomas}. This proves the integrality of BPS states \(N_{S^r,i,j}\), which is the statement of the LMOV conjecture. Alternate proofs of the DT-invariant integrality structure exist in the literature using techniques such as: combinatorics of zonotopical algebras \cite{10.1093/imrn/rnad033} and quiver diagonalisation method \cite{Jankowski_2023}. The numerical DT-invariants play a crucial role in the study of BPS degeneracies associated to quiver gauge theory for a symmetric quiver\footnote{\(s\) plays the role of spin of a BPS particle \cite{Kucharski_2017, Kucharski_2019, Ekholm:2018eee} 
labeled by dimensional vector \(\textbf{d}=(d_1,...,d_m)\).}.
\\ \\
According to \cite{Kucharski_2017,Kucharski_2019}, the motivic generating series assigned to any symmetric quiver takes the following form,
\begin{equation}
    P^{Q_K}(x,q)=\sum_{d_1,..,d_m \geq 0}(-q)^{\sum_{i,j=1}^mC_{ij}d_id_j}\prod_{i=1}^m\frac{x_i^{d_i}}{(q^2;q^2)_{d_i}} \label{PQ}
\end{equation}
where \(C_{ij}\) is the adjacency matrix, and \(d_i\) is the dimensional vector space assigned to each node \(i\), with the total number of nodes of the quiver being \(m\). Moreover, to each node \(i\), a quiver assigns a linear map \(\alpha:i \to j\) to each arrow. Similar to \eqref{pleth}, we can write the quiver generating series as,
\begin{equation}
    P^{Q_K}(x,q)=\prod_{\substack{d_1,..,d_m \geq 0 \\ s \in \mathbb{Z}}}\ \prod_{k \geq 0}(1-x_1^{d_1}...x_m^{d_m}q^{2k+s})^{\Omega^{Q_K}_{d_1...,d_m;s}}
\end{equation}
A geometric interpretation of the DT-invariant coefficients was given \cite{meinhardt2016donaldsonthomasinvariantsversusintersection, Franzen_2018} in terms of the intersection Betti numbers of the moduli space of all semisimple representations of the quiver \(Q_K\) of dimensional vector \(\textbf{d}\), or alternatively as Chow-Betti numbers of moduli space of all simple representations for \(Q_K\) having dimensional vector \(\textbf{d}\).
\\ \\
The main conjecture of the knot-quiver correspondence states that the generating function of the colored HOMFLY-PT polynomial can be written as,
\begin{equation}
    P(x)=\sum_{r=0}^{\infty}\overline{P}_r(a,q)x^r= \sum_{d_1,..,d_m \geq 0}q^{\sum_{i,j=1}^mC_{ij}d_id_j}\frac{\prod_{i=1}^m x_i^{d_i}q_i^{l_id_i}a_i^{a_id_i}(-1)^{t_id_i}}{\prod_{i=1}^m(q^2;q^2)_{d_i}} \label{kqseries} 
\end{equation}
where the quiver matrix \(C_{ij}\) is an \(m \times m\) symmetric matrix, \(l_i,a_i\) and \(t_i\) are fixed integers, encoding the information of the generators of the uncolored HOMFLY-PT homology generators. The term proportional to \(x^r\) arises from a set of dimensional vectors \(\{d_i\}\) with the constraint \(r=d_i+...+d_m\). It is interesting to note that \eqref{PQ} and \eqref{kqseries} are identified upto a special class of variable change,
\begin{equation}
    x_i=xa^{a_i}q^{l_i-1}(-1)^{t_i} \label{kqvariable}
\end{equation}
where the size of the matrix \(C_{ij}\) is equal to the number of generators of the uncolored HOMFLY-PT homology. The integer \(t_i\) encodes the homological degrees of generators of the uncolored HOMFLY-PT homology, which remarkably is also equivalent to the diagonal entries of the matrix \(C_{ij}\) i.e. \(C_{ii}=t_i\) (counts loops of the \(i^{th}\) node). The powers of \(q\) are \(l_i=q_i-t_i\), \(a_i\) is the \(a-\)degrees of generators of the uncolored HOMFLY-PT homology. Moreover, the additional minus sign which comes with the grading \(t_i\) becomes relevant for odd \(t_i\)'s. The KQ-correspondence can further be refined using the quadruply-graded homology \cite{gorsky2013quadruplygradedcoloredhomologyknots}.
\\ \\
We have listed the quantities that match under the KQ-correspondence in Figure \ref{KQ_table}. However, there is another knot operation called “framing" which can be seen on the quiver side. The framing operation by \(f \in \mathbb{Z}\) changes the colored HOMFLY-PT polynomial by a factor, which on the symmetric quiver side translates to,
\begin{equation}
    a^{2fr}q^{fr(r-1)} \label{framing}
\end{equation}
From the perspective of the quiver series, the quadratic piece (in \(r\)) of \(q\) is,
\begin{equation}
    q^{fr^2}=q^{f(\sum_i d_i)^2}=q^{f\sum_{i,j} d_id_j}
\end{equation}
and is equivalent to shifting all the elements of the quiver matrix \(C_{ij}\) by a factor \(f\),
\begin{equation}
    C_{ij} \to C_{ij}+f\begin{bmatrix}
        1&1&...&1 \\
        1&1&...&1 \\
        .&.&...& .& \\
        .&.&...& .& \\
        1&1&...&1 \\
    \end{bmatrix} 
\end{equation}
From the quiver side this operation is equivalent to adding \(f\) loops at each vertex and \(f\) pairs of oppositely-oriented arrows between pairs of all vertices. 
\\ \\
\textbf{Example: Trefoil Knot:} In the case of the trefoil knot \(3_1\), the reduced superpolynomial takes the form \cite{Fuji_2013},
\begin{equation}
    \mathscr{P}^{3_1}_r(a,q,t)=\frac{a^{2r}}{q^{2r}}\sum_{k=0}^r \left[ \begin{array}{l}
  r \\
  k \\
\end{array} \right]
q^{2k(r+1)}t^{2k}\prod_{i=1}^k(1+a^2q^{2(i-2)}t) \label{trefoil}
\end{equation}
where the symbol \(\left[ \begin{array}{l}
  r \\
  k \\
\end{array} \right] = \frac{(q^2;q^2)_r}{(q^2;q^2)_k(q^2;q^2)_{r-k}}\) and the \(q\)-Pochhammer symbol is defined to be \((z;q)_n=\prod_{i=0}^{n-1}(1-zq^i)\). The reduced uncolored \((r=1)\) HOMFLY-PT superpolynomial is given by,
\begin{equation}
    \mathscr{P}^{3_1}_1(a,q,t)=\frac{a^2}{q^2}+a^2q^2t^2+a^4t^3
\end{equation}
Taking the Euler characteristic \((t=-1)\) of \eqref{trefoil} and using the q-binomial identity \((z;q)_k=\sum_i\left[ \begin{array}{l}
  k \\
  i \\
\end{array} \right](-z)^iq^{\frac{i(i-1)}{2}}\), we can massage the right-most product in \eqref{trefoil} to the form,
\begin{equation}
    \left[ \begin{array}{l}
  r \\
  k \\
\end{array} \right]\left(\frac{a^2}{q^2};q^2\right)_k=\sum_{i=0}^k\frac{(q^2;q^2)_r\left(-\frac{a^2}{q^2}\right)^iq^{i(i-1)}}{(q^2;q^2)_{r-k}(q^2;q^2)_i(q^2;q^2)_{k-i}}
\end{equation}
\begin{figure}
    \centering
    \includegraphics[width=0.5\linewidth]{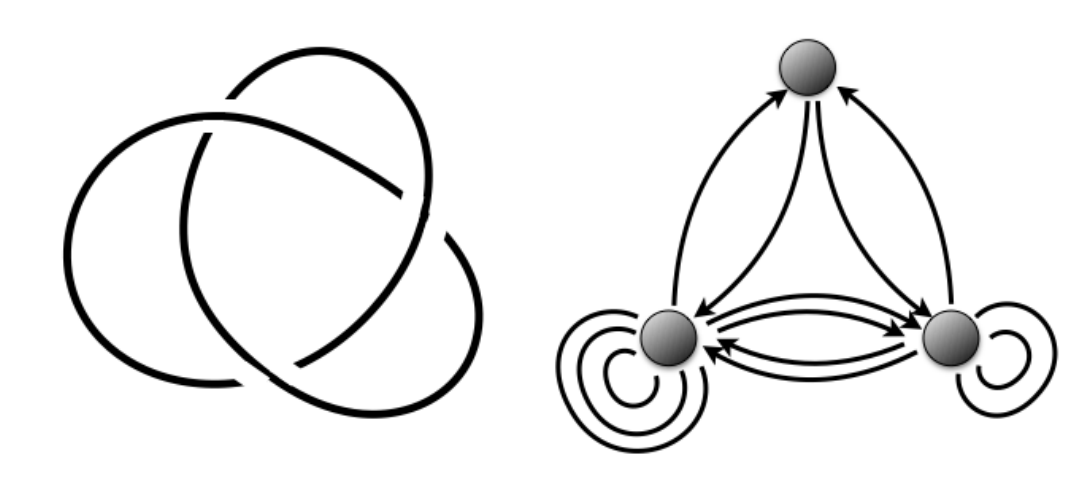}
    \caption{Trefoil knot and its corresponding quiver \cite{Kucharski_2017,Kucharski_2019}.}
    \label{fig:trefoil_quiver}
\end{figure}
Introducing new variables \(r=d_1+d_2+d_3, k+d_2+d_3,i=d_3\) such that \(d_i \geq 0\) and normalizing \(\mathscr{P}^{3_1}_r(a,q,t=-1)\) by \((q^2;q^2)_r\), the HOMFLY-PT generating series \eqref{p(x)} can be rewritten in terms of the quiver series \eqref{kqseries},
\begin{equation}
     P(x)=\sum_{r=0}^{\infty}\frac{P^{3_1}_r(a,q)}{(q^2;q^2)_r}x^r=\sum_{d_1,d_2,d_3 \geq 0}q^{\sum_{i,j=1}^3C_{ij}d_id_j-2d_1-3d_3}\frac{(-1)^{d_3}a_i^{2d_1+2d_2+4d_3} x^{\sum_i d_i}}{(q^2;q^2)_{d_1}(q^2;q^2)_{d_2}(q^2;q^2)_{d_3}}
\end{equation}
where the quiver matrix for the trefoil is
\[
C=\begin{bmatrix}
0 & 1 & 1\\
1 & 2 & 2\\
1 & 2 & 3
\end{bmatrix}
\]
The diagonal elements \(t_i=(0,2,3)\) of the quiver matrix \(C_{ij}\) correspond to the number of loops at the respective nodes of the quiver (homological degree of the knot superpolynomial for \(3_1\)). The power of variable \(q\) having linear terms in \(d_i\) correspond to \(l_i=(-2,0,3)\), and the coefficients \(a_i=(2,2,4)\) agree with the homological degree of the HOMFLY-PT superpolynomial for the trefoil knot. 
\\ \\
Representing the HOMFLY-PT superpolynomial in terms of a quiver \eqref{kqseries}, guarantees that the LMOV invariants can be expressed in terms of the motivic DT invariants \(\Omega^{3_1}_{d_1,d_2,d_3;s}\), proving the LMOV conjecture for the trefoil knot, for all symmetric representations.
\subsection{Unreduced Quivers}
Let's try to briefly understand the role of normalization for the quiver parameters \(a_i, q_i\) and \(C_{ii}\). The duality between the HOMFLY-PT homology generators and the terms in the quiver series only hold when the normalization choice of the unreduced HOMFLY-PT polynomial \(\overline{P}_r(a,q)\) is \((q^2;q^2)_r\) i.e.
\begin{equation}
    \overline{P}_r(a,q)=\frac{P_r(a,q)}{(q^2;q^2)_r} \label{norm1}
\end{equation}
However, choosing to normalize the unreduced polynomial with the full unknot polynomial,
\begin{equation}
    \overline{P}_r(a,q)=a^{-r}q^{r}\frac{(a^2;q^2)_r}{(q^2;q^2)_r}P_r(a,q) 
\end{equation}
produces a twice larger quiver matrix, encoding the information about unreduced HOMFLY-PT homology, and the Poincaré polynomial is obtained by multiplying the reduced superpolynomial by \(a^{-1}q(1+a^2t)\). Suppose the normalization choice \eqref{norm1} leads to the quiver matrix \(C_{ij}\). Introducing an extra factor of \(a^{-r}q^{r}(a^2;q^2)_r\), we use the “splitting" identity 
\begin{equation}
    \frac{(\xi;q^2)_r}{(q^2;q^2)_r}=\sum_{\alpha+\beta=n}(-q)^{\beta^2}\frac{(\xi q^{-1})^{\beta}}{(q^2;q^2)_{\alpha}(q^2;q^2)_{\beta}} \label{split}
\end{equation}
to deal with the additional factor \(q\)-Pochhammer \((a^2;q^2)_r\). Furthermore, we introduce additional variables \(\alpha_i\) and \(\beta_i\) such that \(d_i=\alpha_i+\beta_i\). A detailed analysis has been given in \cite{Kucharski_2019}. The result of this analysis states that the exponent of the \(q-\)variable is as follows:
\begin{equation}
    \sum_{i,j}C_{ij}\beta_i\beta_j+ \sum_{i,j}(C_{ij}+1)\alpha_i\alpha_j+2\sum_{i\leq j}C_{ij}\alpha_i\beta_j+2\sum_{i > j}(C_{ij}+1)\alpha_i\beta_j \label{largerC}
\end{equation}
This expression encodes a quiver of twice the size of the original quiver \(C_{ij}\). The new quiver \(\widetilde{C}_{i'j'}\) decomposes into a \(C_{ij}\) like term: \(\left(\sum_{i,j}C_{ij}\beta_i\beta_j\right)\) and another term which also has the original quiver structure but shifted by framing 1 \(\left(\sum_{i,j}(C_{ij}+1)\alpha_i\alpha_j\right)\). These two subquivers are connected via arrows, whose structure is given in the last two terms of \eqref{largerC}.
\\ \\
\textbf{Case study (Unknot):} The unreduced colored HOMFLY-PT polynomial for the unknot is of the form,
\begin{equation}
     \overline{P}_r(a,q)=a^{-r}q^{r}\frac{(a^2;q^2)_r}{(q^2;q^2)_r} \label{funknot}
\end{equation}
More generally, we can also study the framed unknot, whose connection to the LMOV invariants was studied in \cite{Kucharski_2016, luo2016integralitystructurestopologicalstrings, zhu2018topologicalstringsquivervarieties}. Considering the full unknot invariants \eqref{funknot} with \(r=d_1+d_2\) and using the splitting formula \eqref{split}, the knot generating function can be written as follows,
\begin{equation}
    P(x)=\sum_{r=0}^{\infty}x^ra^{-r}q^{r}\frac{(a^2;q^2)_r}{(q^2;q^2)_r}=\sum_{d_1,d_2 \geq 0}x^{d_1+d_2}\frac{(-q)^{d_2^2}(a^2q^{-1})^{d_2}}{(q^2;q^2)_{d_1}(q^2;q^2)_{d_2}} \label{funknotpx}
\end{equation}
Thus the new quiver \(\widetilde{C}_{i'j'}\) and the generators \(x\) are given by,
\[
\widetilde{C}_{i'j'}=\begin{bmatrix}
0 & 0 \\
0 & 1 \\
\end{bmatrix},
\quad
x= \begin{bmatrix}
x \\
a^2q^{-1}x
\end{bmatrix}
\]
The larger quiver \(\widetilde{C}_{i'j'}\) consist of two disconnected nodes, with one of nodes (labeled by \(d_2\) in this example) having a single loop. Moreover writing \eqref{funknotpx} in terms of quantum dilogarithms gives two non-zero LMOV or equivalently motivic DT invariants, which is true for an unknot \cite{Ooguri_2000,Kucharski_2016}.
\\ \\
Including the framing dependence \eqref{framing} in \eqref{funknotpx} again leads to a quiver with additional loops. However, unlike the unframed unknot \eqref{funknotpx} the result will not factorize in finite quantum dilogarithms, thus encoding an infinite number of LMOV invariants.
\\ \\
As our final comment of this section we will briefly be introducing the quiver operation of (un)linking and involution \cite{Kucharski:2023jds}. Contrary to the framing operation, (un)linking and involution preserves the quiver generating series and thus the motivic Donaldson-Thomas invariants. Due to such operators, knots and quivers are not bijective. Furthermore, using (un)linking, one can introduce a “permutohedra" like structure for quivers,\footnote{A neat introduction to knot-quiver correspondence and various operators related to quivers can be found in the lectures of Dmitry Noshchenko delivered during the Simons semester “Knots, Homologies and Physics" at Warsaw, 2024.} whose edges represent quivers which differ by a single transposition and can be unlinked once to yield the same quiver \cite{Jankowski:2021flt}.
\section{Geometry of the knots-quivers duality} \label{Geometry_KQ}
In this section we will view the knot-quiver correspondence from an alternative perspective, relating them directly to geometric and physical objects. We will try give a brief overview of how the relationship between knots and quivers can be seen from open topological string theory via brane dynamics as well as from the geometry of holomorphic curves attached to branes.
\begin{figure}[h]
    \centering
    \includegraphics[width=0.7\linewidth]{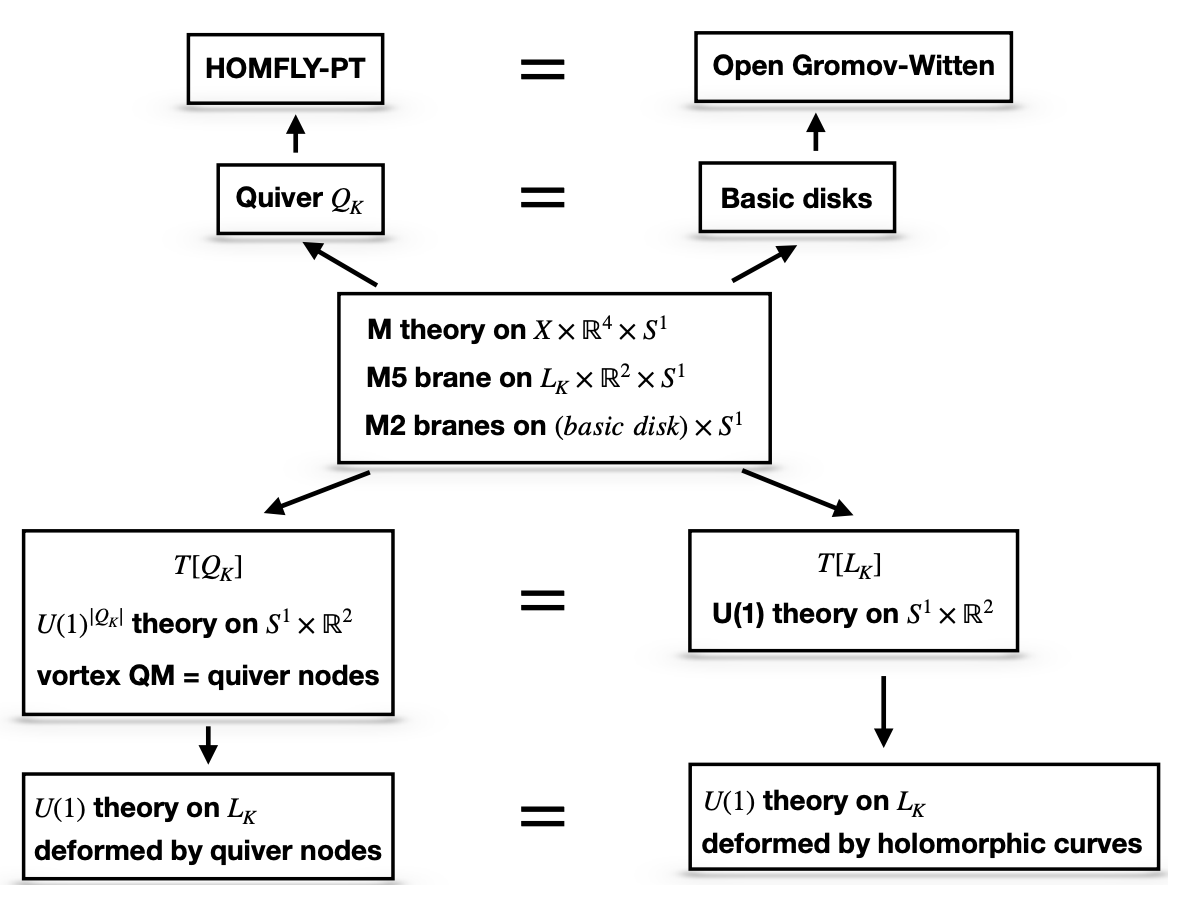}
    \caption{An equivalence between the respective geometric and physics side of quivers \cite{ekholm2020physicsgeometryknotsquiverscorrespondence}.}
    \label{fig:geometry-KQ}
\end{figure}
\subsection{Physical picture}
We are aware that under the large \(N\) transition, the colored HOMFLY-PT polynomial can be described using the Gromov-Witten invariants of the knot conormal \(\mathscr{L}_K\) in the resolved conifold \(X\) \cite{Ooguri_2000, Vafa_2014, ekholm2025skeinsbranes}. Through the lens of M-theory, the generating series of the HOMFLY-PT polynomial counts M2-branes wrapping the holomorphic disks and ending on the M5-brane wrapping \(\mathscr{L}_K\). The full M-theory background is \(X \times S^1 \times \mathbb{R}^4\), of which the M5-brane wrap \(\mathscr{L}_K \times S^1 \times \mathbb{R}^2\), while the M2-branes wrap the product of the holomorphic curve and one-cycle \(S^1\). Physically this describes a \(U(1)\) Chern-Simons theory on the Lagrangian conormal deformed by the embedded holomorphic disks. However, one can also view this theory as a 3d \(\mathcal{N}=2\) theory \(T[\mathscr{L}_K]\) on \(S^1\times \mathbb{R}^2\) \cite{Dimofte_2011, Terashima_2011,Yagi_2013,Lee_2013, cordova2013complexchernsimonsm5branessquashed}.
\\ \\
To a theory \(T[\mathscr{L}_K]\), there exists a dual description which features a \(U(1)\) group and a single chiral field assigned to each node of the quiver. Interaction between the nodes are controlled by the Chern-Simons and the Fayet-Ilioupoulos (FI) couplings. The FI term couples the fugacities of the global and its topological dual symmetry in \(T[Q_K]\). Moreover, in \(T[Q_K]\) theory the wrapped M2-branes with endings on M5-branes count the BPS vortices, which admits a quantum mechanics description \cite{Hwang_2017}. The quiver description of the BPS spectrum can be understood by identifying the M2-branes with quiver nodes and intersections between M2-branes as the arrows connecting the nodes. 
\subsection{Geometric picture}
We now introduce the geometric description of the quivers. The quiver nodes are in one to one correspondence to the embedded holomorphic disks with boundary on \(\mathscr{L}_K\) with a framing along the knot. The basic disks wrapped around \(\mathscr{L}_K\) deform the symplectic structure of the Lagrangian \(T^*\mathscr{L}_K\). Each node in the dual theory \(T[Q_K]\), consists of a \(U(1)\) gauge theory and a fundamental chiral, corresponding to a \(U(1)\) CS-theory deformed by all the generalized holomorphic curves ending on \(\mathscr{L}_K\).
\\ \\
The winding number \(a_i\) of the \(i^{th}\) basic disk around \(\mathbb{P}^1\) of the resolved conifold \(X\), \(q_i=e^{g_s/2}\) is related to the self linking of the \(i^{th}\) basic disk, \(l_i=q_i-C_{ii}\) counts the intersection number of the basic disks with the 4-chain \(C_K\) with boundary \(\partial C_K=2\mathscr{L}_K\) and \(t_i=C_{ii}\) is the wrapping number of the \(i^{th}\) basic disk around the \(S^2-\)cycle connecting \(\mathscr{L}_K\). The 4-chain \(C_K\) has its origins in knot contact homology \cite{ekholm2020highergenusknotcontact} and a physical interpretation of this 4-chain has been illustrated in \cite{ekholm2020physicsgeometryknotsquiverscorrespondence}. Precisely, this can be summarized in the following way: The refined Gromov-Witten partition function \(\Psi_K(x,a,t,q)\) counts all the (disconnected) generalized holomorphic curves ending on \(\mathscr{L}_K\) and is expressed via a generating series,
\begin{equation}
    \Psi_K(x,a,t=-1,q)=\sum_{r \geq 0}P^K_r(a,q)x^r
\end{equation}
where \(P^K_r(a,q)\) is the HOMFLY-PT polynomial colored with the symmetric representations.   
\\ \\
Let there be an embedded knot \(K \subset S^3\) and a Lagrangian conormal \(\mathscr{L}_K \subset X\) which is shifted off the zero-section of \(T^*S^3\) in the resolved conifold \(X\). For a finite number of embedded disks \(D_i,...,D_m\) (as many as the number of nodes in the symmetric quiver \(Q_K\)) with endpoints on \(\mathscr{L}_K\), the \(D_i\) are the nodes and the arrows \(C_{ij}\) are the (self)linking numbers \(\text{lk}(\partial D_i,\partial D_j)\) counted with respect to a reference framing of \(\partial \mathscr{L}_K\). The function \(x_i\) is the quiver generating variable of \(D_i\), \(a_i\) is the homology class of \(D_i \in H_2(X)\), \(n_i\) denotes the homology class of \(\partial D_i \in H_1(X)\) and \(l_i=q_i-C_{ii}\) denotes the intersection \(D_i.C_K\). Using this prescription we get,
\begin{equation}
    \Psi_K(x,a,t,q)=P^{Q_K}(x_1,...,x_m)|_{x_i=x^{n_i}a^{a_i}q^{l_i}(-t)^{C_{ii}}}
\end{equation}
where \(P^{Q_K}\) is the quiver generating series \eqref{PQ}. A similar relationship also connects the LMOV invariants and the motivic DT invariants,
\begin{equation}
    \Omega^{Q_K}(x_i,q)=|_{x_i=x^{n_i}a^{a_i}q^{l_i}(-t)^{C_{ii}}}=N^K(x,a,q)
\end{equation}
\subsection{Holomorphic disks and Superpotential}
Given a Lagrangian conormal \(\mathscr{L}_K \subset X\) shifted off the zero-section of \(T^*S^3\) for an embedded knot \(K \subset S^3\), it was shown that the corresponding HOMFLY-PT polynomial for the knot \(K\) is identified with the open Gromov-Witten invariants of \(\mathscr{L}_K \subset X\), i.e the wavefunction \(\Psi_K(\xi)\),
\begin{equation}
    \Psi_K(\xi)=\text{exp}\left(\sum_{n \geq 1}C_n(e^t, g_s)e^{n\xi} \right) \label{gwwf}
\end{equation}
where \(C_n(e^t, g_s)\) is a polynomial in \(e^t\) that counts (connected) generalized holomorphic curves with boundary on \(\mathscr{L}_K\) \cite{ekholm2020highergenusknotcontact} and \(t\) is in the homology class \(H_1(X)\). In the semi-classical limit \(g_s \to 0\), the wavefunction \eqref{gwwf}, one gets is the Gromov-Witten disk potential computed via the Abel-Jacobi map on the mirror curve \cite{aganagic2000mirrorsymmetrydbranescounting, Aganagic_2002}
\begin{equation}
    \Psi_K(\xi) \sim \text{exp}\left(\frac{1}{g_s}W_K(\xi)+...  \right) \label{GWdiskpot}
\end{equation}
In \cite{Vafa_2014, ekholm2020highergenusknotcontact} a different approach to calculate the Gromov-Witten wavefunction was identified using Knot Contact Homology. Using this construction, the mirror curve matched with the augmentation curve Aug\(_K\) of Knot Contact Homology. 
\\ \\
We now move onto the extraction of the twisted superpotential and its Legendre transformation. In the theory \(T[\mathscr{L}_K]\) and its weak coupling counterpart \(T_0[\mathscr{L}_K]\), the moduli space of supersymmetric vacua coincides with the moduli space of flat connections on \(\mathscr{L}_K\). The twisted superpotential of \(T_0[\mathscr{L}_K]\) is encoded in the large color and semi-classical limit of the colored superpolynomial,
\begin{equation}
   \mathscr{P}^K(a,q,t) 
    \xrightarrow[{r \to \infty}]{\hbar \to 0}  \int \prod_i \frac{dz_i}{z_i}\text{exp}\frac{1}{2\hbar}\left(\widetilde{W}_{T_0[\mathscr{L}_K]}(z_i,a,t,y)\ +...\right) \label{T_0[L_K]}
\end{equation}
Under this limit, we have \(y \to q^{2r} \). The twisted superpotential \(\widetilde{W}_{T_0[\mathscr{L}_K]}\) is a function of several fugacities: 
\begin{itemize}
    \item \(y\) is associated to the global symmetry \(U(1)_M\) and is identified with the meridian holonomy of flat connections on \(\mathscr{L}_K\) or equivalently the torus \(T^2\) at infinity.
    \item \(a\) is associated to the global symmetry \(U(1)_Q\) arising due to the 2-cycle \(S^2\) of \(\mathbb{CP}^1\) of the resolved conifold \(X\).
    \item \(-t\) is associated to the rotational symmetry \(U(1)_F\) of the normal bundle \(\mathbb{R}^2 \subset \mathbb{R}^4\).
    \item \(z_i\) is associated to the abelian gauge symmetry \(U(1) \times ... \times U(1)\). These fugacities are identified by \(z_i=q^{2k_i}\).
\end{itemize}
The twisted superpotential includes chiral and CS-coupling terms denoted by,
\begin{align}
    & \text{Li}_2(a^{n_Q}(-t)^{n_F}y^{n_M}z_i^{n_i}) \leftrightarrow \text{(chiral field)} \nonumber \\ & \frac{\kappa_{ij}}{2}\log{\theta_i}\log{\theta_j} \leftrightarrow \text{(CS coupling)}
\end{align}
The linear dilogarithm term is understood as a one-loop contribution of a chiral field charged under the various \(U(1)\) fugacities, while the quadratic dilogarithm term denotes the CS-coupling term between different fugacities \(\theta_i\). 
\\ \\
The effective superpotential \(\widetilde{W}^{\text{eff}}_{T_0[\mathscr{L}_K]}(a,t,y)\) is extracted by gauging out the abelian gauge symmetry, which is identical to taking the saddle point approximation with respect to \(z_i\). In \cite{Fuji_2013}, it was shown that the moduli space of supersymmetric vacua of the theory \(T_0[\mathscr{L}_K]\) coincides with the graph of the super \(A-\)polynomial,
\begin{equation}
    \frac{\partial \widetilde{W}^{\text{eff}}_{T_0[\mathscr{L}_K]}(a,t,y)}{\partial \log{y}}=-\log{x} \leftrightarrow \mathcal{A}^K(x,y,a,t)=0
\end{equation}
Let's return to the theory of our interest \(T[\mathscr{L}_K]\), and take a variant of the double-scaling limit \eqref{T_0[L_K]},
\begin{equation}
    \mathscr{P}^K(a,q,t) 
    \xrightarrow[\substack{q^{2r} \to y \\ q^{2k_i} \to z_i}]{\hbar \to 0}  \int \frac{dy}{y}\int \prod_i \frac{dz_i}{z_i}\text{exp}\frac{1}{2\hbar}\left(\widetilde{W}_{T[\mathscr{L}_K]}(z_i,a,t,y)\ +...\right) \label{T[L_K]}
\end{equation}
where,
\begin{equation}
    \widetilde{W}_{T[\mathscr{L}_K]}(z_i,a,t,y)=\widetilde{W}_{T_0[\mathscr{L}_K]}(z_i,a,t,y)+\log{x}\log{y}
\end{equation}
The FI coupling exists between the \(U(1)_M\) global symmetry associated to the fugacity \(y\) with its dual topological symmetry \(U(1)_L\) having fugacity \(x\). The variable \(x\) is the source term for the BPS vortices. Combining the saddle-point approximation first with respect to the fugacity \(z_i\) and then \(y_i\), 
\begin{equation}
    \frac{\partial \widetilde{W}^{\text{eff}}_{T[\mathscr{L}_K]}(a,t,y)}{\partial \log{y}}=0 \leftrightarrow \mathcal{A}^K(x,y,a,t)=0
\end{equation}
where 
\begin{equation}
    \widetilde{W}^{\text{eff}}_{T[\mathscr{L}_K]}(x,a,t,y)=\widetilde{W}^{\text{eff}}_{T_0[\mathscr{L}_K]}(a,t,y)+\log{x}\log{y}
\end{equation}
\\ \\
To illustrate how the Gromov-Witten disk potential is related to the twisted superpotential of the 3d \(\mathcal{N}=2\) theory \(T[\mathscr{L}_K]\), recall that the semi-classical approximation of the open-topological string wavefunction \(\Psi_K\) gives the GW disk potential \eqref{GWdiskpot},
\begin{equation}
     P_K(x,a,q) \stackrel{\hbar \to 0}{\sim} \exp\left(\frac{1}{2\hbar}W_K(a,x)+...\right)
\end{equation}
Whereas, taking the limit \(\hbar \to 0\) of the HOMFLY-PT generating series as defined in \eqref{T[L_K]} along with the saddle-point approximations and gauging out both the fugacities \(z_i\) and \(y_i\) gives us the Legendre transform,
\begin{equation}
  W_K(a,x)=\left(\widetilde{W}^{\text{eff}}_{T_0[\mathscr{L}_K]}(a,t,y)+\log{x}\log{y}\right)|_{y=y^*(x)} \label{legendreWK}
\end{equation}
where \(y^*(x)\) is the saddle point. There is a quantum analog of the statement that \(W_K\) is the Legendre transform of \(\widetilde{W}^{\text{eff}}_{T[\mathscr{L}_K]}\), where the full quantum effective action of the gauge theory is the Fourier transform of the open-topological string wavefunction \(\Psi_K\).
\\ \\
Solving the zero locus equation for the \(\mathcal{A}^K(x,y,a)=0\) for \(y\) gives a function involving the classical LMOV invariants \eqref{classlmov} \(b^K_{r,i}=\sum_jN^K_{r,i,j}\) \cite{Garoufalidis_2016},
\begin{equation}
    y^*(a,x)= \lim_{\hbar \to 0} \frac{P^K(q^2x,a,q)}{P^K(x,a,q)}
\end{equation}
Integrating this function gives the Gromov-Witten disk potential,
\begin{equation}
    W_K(a,x)= \int d\log{x}\log{y^*}
\end{equation}
while solving for the function \(x^{-1}\) and integrating with respect to it gives the twisted superpotential of \(T_0[\mathscr{L}_K]\)
\begin{equation}
    \widetilde{W}^{\text{eff}}_{T_0[\mathscr{L}_K]}=\int d\log{y}\ (\log{x^{-1}})^*
\end{equation}
\subsection{Quiver theory \(T[Q_K]\)}
The essence of the knot-quiver correspondence involves expressing the knot polynomials in terms of the representation theory of quivers. As reviewed in Section \ref{Knots&Quivers}, we have 
\begin{equation}
    \mathscr{P}(a,q,t,x)=P^{Q_K}(x_i,q)|_{x_i=xa^{a_i}q^{l_i}(-1)^{C_{ii}}}
\end{equation}
where both sides of the equation can be written as \eqref{PQ}. If we apply \eqref{T[L_K]} directly to \(P^{Q_K}\), we obtain the dual 3d \(\mathcal{N}=2\) theory \(T[Q_K]\) given by,
\begin{equation}
    P^{Q_K}(x_i,q) 
    \xrightarrow[\substack{q^{2d_i} \to y_i}]{\hbar \to 0} \int \prod_{i \in Q_{K_0}} \frac{dy_i}{y_i}\text{exp}\frac{1}{2\hbar}\left(\widetilde{W}_{T[Q_K]}(x_i,y_i)\ +...\right) \label{T[Q_K]}
\end{equation}
where the quiver superpotential is,
\begin{equation}
    \widetilde{W}_{T[Q_K]}(x_i,y_i)=\sum_i\text{Li}_2(y_i)+\log{((-1)^{C_{ii}}x_i)}\log{y_i}+\sum_{i,j}\frac{C_{ij}}{2}\log{y_i}\log{y_j}
\end{equation}
\(\widetilde{W}_{T[Q_K]}\) encodes a finite set of sources (one for each node \(Q_{K_0}\) together with a finite set of couplings (\(C_{ij}\) are the arrows \(i \to j\).) The quiver dictionary reads as follows,
\begin{enumerate}
    \item Gauge group: \(U(1)^{(1)} \times ... \times U(1)^{(m)}\)
    \item Matter: one chiral field \(\phi_i\) at each node, charged \(\delta_{ij}\) under \(U(1)^{(i)}\)
    \item mixed CS coupling: \(\kappa^{\text{eff}}_{ij}=C_{ij}\)
    \item FI couplings: \(\log{((-1)^{C_{ii}}x_i)}\)
\end{enumerate}
The semi-classical limit of the quiver generating series gives the quiver disk potential, described as
\begin{equation}
    W_{Q_K}(x_i)=\lim_{\hbar \to 0}2\hbar \log{\left(P^{Q_K}(x_i,q)\right)}
\end{equation}
Once again, under the Legendre transform, it is identified with the quiver twisted superpotential,
\begin{equation}
     W_{Q_K}(x_i)=\widetilde{W}_{T[Q_K]}(x_i,y_i)|_{y_i=y_i^*(x_i)}
\end{equation}
The saddle-point \(y_i^*(x_i)\) defines the quiver \(A-\)polynomial \(\mathcal{A}_i^{Q_K}(x_i,y_i)\) analogous to \cite{Panfil_2018, Panfil_2019},
\begin{equation}
    \frac{\partial \widetilde{W}^{\text{eff}}_{T[Q_K]}(a,t,y)}{\partial \log{y_i}}=0 \leftrightarrow \mathcal{A}_i^{Q_K}(x_i,y_i)=0 \leftrightarrow \frac{\partial W_{Q_K}(x_i)}{\partial \log{x_i}}=\log{y_i}
\end{equation}
Using the same analysis as before, we arrive at 
\begin{equation}
    \mathcal{A}_i^{Q_K}(x_i,y_i)=1-y_i-x_i(-y_i)^{C_{ii}}\prod_{i \neq j}y^{C_{ij}}_j
\end{equation}
Moreover, in the semi-classical limit of the KQ variables, we get a relation between the Gromov-Witten and quiver disk potential,
\begin{equation}
    W_{Q_K}(x_i)|_{x_i=xa^{a_i}}=W_K(a,x) \label{WW}
\end{equation}
From \eqref{WW} and \eqref{legendreWK}, we are able to establish an interesting relation between the \(A-\)polynomial for a knot \(K\) and it's quiver counterpart \(\mathcal{A}_i^{Q_K}\),
\begin{equation}
    y^*(x)=\prod_iy^*_i(x_i)|_{x_i=xa^{a_i}} \label{meridian}
\end{equation}
where the term \(y^*(x)\) solves for the zero locus equation \(\mathcal{A}^K(x,y)=0\) and \(y^*_i(x_i)\) solves \(\mathcal{A}_i^{Q_K}(x_i,y_i)\). Geometrically,  \(y^*_i(x_i)\) is the meridian holonomy for the flat \(U(1)\) connections on \(\mathscr{L}_K\) (around each source \(x_i\)) on the tubular neighborhoods of basic disks and their sum adds up to the meridian holonomy \(y^*(x)\) of the torus \(T^2\) at infinity.
\\ \\
We conclude this section by giving a geometric interpretation of the quiver data i.e. the variables \(a_i, q_i,t_i\) of the KQ correspondence. The topological data that the holomorphic disks encodes is the homology class of the disk boundary \(\partial \Sigma \in H_1(\mathscr{L}_K)\) and the wrapping number of the disks around the two-cycles in \(H_2(X)\). The knot conormal is topologically equivalent to \(S^1 \times \mathbb{R}^2\), thus basic disks wrap the one-cycle exactly once, justifying why in the KQ change of variables \(x_i \sim x\). It is worth noting that for knots which don't satisfy the exponential growth property \eqref{expo_sp}, the KQ correspondence still holds, but with a modified KQ variable change \(x_i \sim x^{n_i}\)  for certain nodes of the quiver. This highlights that in the most general cases, we also need to consider the contributions from multiply-wrapped basic disks. One obvious example, where such modification was first used is the knot \(9_{42}\) \cite{ekholm2020physicsgeometryknotsquiverscorrespondence, ekholm2022knothomologiesgeneralizedquiver}. Next, the variable \(a_i\) counts the wrapping number of the \(i^{th}\) basic disk around \(\mathbb{CP}^1\). The KQ variable \(t_i=C_{ii}\) corresponds to the winding number of the \(i^{th}\) basic disk around around the \(S^2\) connecting \(\mathscr{L}_K\). Finally the variable \(q_i\) is responsible for counting the self-linking of the \(i^{th}\) basic disk. The powers of \(q\) always appear in the combination \((q_i-C_{ii})\), which equals the intersection number of the disk with the 4-chain.
\subsection{Quivers and Vortices}
In this section we will show that the BPS vortices partition function of \(T[Q_K]\) coincides with the generating function of the motivic Donaldson-Thomas invariants of a quiver \(Q_K\). 
\\ \\ 
For a simple quiver having only diagonal entries \(C_{ij}=\kappa^{\text{eff}}_i\delta_{ij}\) (disconnected nodes). The effective theory \(T[Q_K]\) consists of \(m-\)copies of \(U(1)_{\kappa_i}\) gauge theory with a single chiral field having charge \(+1\) at each node. The vortex partition function \(\mathcal{Z}^{\text{vortex}}\) can be written as \cite{Beem_2014, hwang2013factorization3dsuperconformalindex, Shadchin_2007}
\begin{equation}
   \mathcal{Z}^{\text{vortex}}=\prod_{j=1}^m\sum_{d_j \geq 0}\left(2(-e^{\mu})^{\kappa^{\text{eff}}_j}(e^{\gamma-\mu})^{1/2}z_j\right)^{d_j}(-e^{\gamma})^{\kappa^{\text{eff}}_jd^2_j}\prod_{s=1}^{d_j}\frac{1}{1-e^{2\gamma}} \label{vortex}
\end{equation}
where \(\gamma, \mu\) are the parameters for the global rotational symmetries \(U(1)_{\gamma} \times U(1)_{\mu}\) for the tangent and the normal bundle to the \(\mathbb{R}^2 \subset \mathbb{R}^4\) respectively and \(d_j\) represents vorticity. Adopting the 3d-3d dictionary, one identifies of variables \(e^{\gamma}=q\) and \(z_i=x_i\), upto some normalization factors, we exactly match the vortex partition function to the quiver partition function,
\begin{equation}
    \mathcal{Z}^{\text{vortex}}=P^{Q_K}
\end{equation}
For a more general quiver (with connected nodes), one needs to construct the vortex partition function in terms of the holomorphic blocks \cite{Beem_2014}. Using this method, the authors in \cite{ekholm2020physicsgeometryknotsquiverscorrespondence} were able to match the vortex partition function to the quiver partition function even for a general quiver \(Q_K\). This establishes the fact that the vortex partition function arises from the vortex quantum mechanics of \(T[Q_K]\). 
\\ \\
The work \cite{Hwang_2017} sheds light on the fact that vortex quantum mechanics of a class of 3d \(\mathcal{N}=2\) theories admit a quiver description. The M-theory description of \(T[\mathscr{L}_K]\) connects basic disks to BPS vortices \cite{Ooguri_2000}, indicating that we must have a dual \(T[Q_K]\) theory. According to \cite{Hwang_2017}, the vortex partition function also admits a quiver description, implying that there must be “fundamental vortices" that generate the entire BPS spectrum. However, using the quiver geometric picture, we can also identity the quiver nodes as basic disks and the spectrum of disks create the boundstates in \(T[Q_K]\). Both these descriptions although different, generate the entire BPS spectrum. In \cite{ekholm2020physicsgeometryknotsquiverscorrespondence}, this relation was extensively studied for the unknot and was shown to be equal to each other. 
\subsection{Working out the Unknot}
From the KQ correspondence \cite{Kucharski_2017, Kucharski_2019}, the unreduced unknot quiver contains two nodes, with a loop on one of the nodes
\begin{figure}[H]
    \centering
    \includegraphics[width=0.4\linewidth]{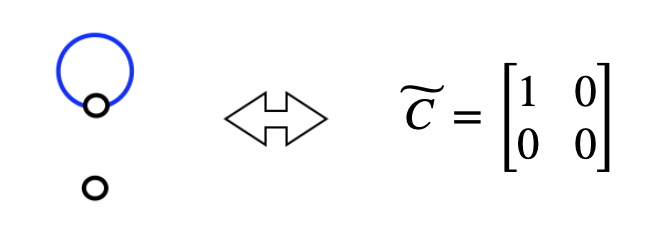}
    \label{fig:Quiver+matrix}
\end{figure}
The quiver generating series for the unreduced unknot \(Q_{0_1}\) is
\begin{equation}
    P^{Q_{0_1}}=\sum_{d_1,d_2 \geq 0}(-q)^{d_1^2}\frac{x_1^{d_1}}{(q^2;q^2)_{d_1}} \frac{x_2^{d_2}}{(q^2;q^2)_{d_2}}\label{funknotpx}
\end{equation}
In terms of the Plethystic exponential,
\begin{equation}
    P^{Q_{0_1}}=\text{Exp}\left(\frac{\Omega^{Q_{0_1}}(x,q)}{1-q^2}\right)
\end{equation}
it gives us two non-zero motivic DT invariants,
\begin{equation}
    \Omega^{Q_{0_1}}_{(1,0),1}=1 \;\ , \;\ \Omega^{Q_{0_1}}_{(0,1),0}=1
\end{equation}
while the numerical DT invariants are as follows,
\begin{equation}
    \Omega^{Q_{0_1}}_{(1,0)}=-1 \;\ , \;\ \Omega^{Q_{0_1}}_{(0,1)}=1
\end{equation}
These values correspond to (uncolored representation) \(|d|=1\), so the BPS states are generated by singly-wrapped basic disks. Using the KQ variable change \(x_1=a^2q^{-1}x,\ x_2=x \), relates the LMOV to the DT invariants,
\begin{align}
    \Omega^{Q_{0_1}}(x,q)|_{x_1=a^2q^{-1}x,\ x_2=x}=&N^{0_1}(x,a,q) \nonumber \\ x-a^2x=&\sum_{r,i,j}N^{0_1}_{r,i,j}x^ra^i,q^j
\end{align}
giving 
\begin{equation}
    N^{0_1}_{1,0,0}=1 \;\ , \;\ N^{0_1}_{1,2,0}=-1
\end{equation} 
while the classical LMOV invariants are,
\begin{equation}
    b^{0_1}_{1,0}=1 \;\ , \;\ b^{0_1}_{1,2}=-1
\end{equation}
Simultaneously, we connects the quiver generating series to the (unreduced) HOMFLY-PT generating series,
\begin{align}
P^{Q_{0_1}}(x_i,q)|_{x_1=a^2q^{-1}x,\ x_2=x} &= P^{{0_1}}(x,a,q) \nonumber \\  &=\sum_{r=0}^{\infty}\frac{(a^2;q^2)_r}{(q^2;q^2)_r}x^r
\end{align}

The same analysis holds for the KQ variables refined by \(t\) or for the doubly refined KQ variables associated to the quadruply graded homology \cite{gorsky2013quadruplygradedcoloredhomologyknots}. 
\\ \\
We now discuss the 3d \(\mathcal{N}=2\) theory for the unknot. We start off with the weak coupling limit theory \(T_0[L_{0_1}]\). Taking the large-color, semi-classical limit \eqref{T_0[L_K]} we get,
\begin{equation}
   \mathscr{P}^{0_1}(a,q,t) 
    \xrightarrow[\substack{q^{2r} \to y}]{\hbar \to 0}  \left[\text{exp}\frac{1}{2\hbar}\left(\widetilde{W}_{T_0[L_{0_1}]}(a,t,y)\ +...\right)\right]
\end{equation}
where 
\begin{equation}
    \widetilde{W}_{T_0[L_{0_1}]}(a,t,y)=\text{Li}_2(y)-\text{Li}_2(-a^2ty)+\text{Li}_2(-a^2t)
\end{equation}
Taking the saddle point with respect to \(y\) gives the zero locus of the super \(A-\)polynomial
\begin{equation}
    \log{x^{-1}}= \frac{\partial \widetilde{W}_{T_0[L_{0_1}]}(a,t,y)}{\partial \log{y}}=-\log{(1-y)}+\log{(1+a^2ty)}
\end{equation}
giving,
\begin{equation}
    \mathcal{A}^{0_1}(x,y,a,t)=1-x-y-a^2txy=0 \label{superA}
\end{equation}
Now we shift our attention to the theory \(T[L_{0_1}]\) \eqref{T[L_K]} gives,
\begin{equation}
    \mathscr{P}^{0_1}(a,q,t) 
    \xrightarrow[\substack{q^{2r} \to y}]{\hbar \to 0}  \int \frac{dy}{y}\left[\text{exp}\frac{1}{2\hbar}\left(\widetilde{W}_{T[L_{0_1}]}(a,t,y)\ +...\right)\right]
\end{equation}
where the twisted superpotential  
\begin{equation}
    \widetilde{W}_{T[L_{0_1}]}(a,t,y)=\text{Li}_2(y)-\text{Li}_2(-a^2ty)+\text{Li}_2(-a^2t)+\log{x}\log{y}
\end{equation}
corresponds to the \(U(1)\) gauge theory having one fundamental and anti-fundamental chiral field. The saddle-point approximation gives,
\begin{equation}
    0=\frac{\partial \widetilde{W}_{T[L_{0_1}]}(a,t,y)}{\partial \log{y}}=-\log{(1-y)}+\log{(1+a^2ty)}+\log{x}
\end{equation}
Solving for \(y\) gives the saddle point \(y^*(x)\). Taking the Legendre transform of \(\widetilde{W}_{T[L_{0_1}]}\) with the specialization \((t=-1)\)  gives the Gromov-Witten disk potential,
\begin{equation}
    W_{0_1}(a,x)=\text{Li}_2(a^2x)-\text{Li}_2(x)
\end{equation}
Moving forward, we shift to the quiver side. Taking the limit \eqref{T[Q_K]} gives,
\begin{equation}
     P^{Q_{0_1}}(x_i,q) 
    \xrightarrow[\substack{q^{2d_i} \to y_i}]{\hbar \to 0}  \int \frac{dy_1}{y_1}\frac{dy_2}{y_2}\left[\text{exp}\frac{1}{2\hbar}\left(\widetilde{W}_{T[Q_{0_1}]}(x_i,y_i)\ +...\right)\right]
\end{equation}
where the quiver potential is 
\begin{equation}
    \widetilde{W}_{T[Q_{0_1}]}=\text{Li}_2(y_1)+\text{Li}_2(y_2)+ \log{(-x_1)}\log{y_1}+\log{x_2}\log{y_2}+\frac{1}{2}\log{y_1}\log{y_1}
\end{equation}
The theory \(T[Q_{0_1}]\) is a \(U(1)^{(1)} \times U(1)^{(2)}\) gauge theory with one chiral field for each node of the quiver. The saddle-points for the quiver superpotential are,
\begin{align}
    &0=\frac{\partial \widetilde{W}_{T[Q_{0_1}]}}{\partial \log{y_1}}=\log{(-x_1)}+\log{y_1}-\log{(1-y_1)} \nonumber \\ &0=\frac{\partial \widetilde{W}_{T[Q_{0_1}]}}{\partial \log{y_2}}=\log{(x_2)}+\log{y_2}-\log{(1-y_2)}
\end{align}
giving the quiver \(A-\)polynomial, 
\begin{align}
    &\mathcal{A}_1^{Q_{0_1}}=1-y_1+x_1y_1 \nonumber \\&\mathcal{A}_2^{Q_{0_1}}=1-y_2+x_2
\end{align}
Notice that \(\mathcal{A}_1^{Q_{0_1}}\) and \(\mathcal{A}_2^{Q_{0_1}}\) are decoupled, which comes from the fact that there are no quiver arrows between the two nodes of the unknot. 
\\ \\
One can also obtain the \(A-\)polynomial from the quiver \(A-\)polynomial. Solving for the zero locus of \(\mathcal{A}^{Q_{0_1}}\), we get 
\begin{equation}
    y_1^*(x_1)=\frac{1}{1-x_1} \;\ , \;\ y_2^*(x_2)=1-x_2
\end{equation}
Using the property \eqref{meridian} yields,
\begin{equation}
    y^*(x)=y_1^*(x_1)y_2^*(x_2)|_{x_1=a^2x,\ x_2=x}=\frac{1-x}{1-a^2x}
\end{equation}
This is the solution to the super \(A-\)polynomial \eqref{superA} for \((t=-1)\). 
The final task will be to relate the quiver disk and the Gromov-Witten disk potential for the unknot in the semi-classical limit \((q \to 1)\),
\begin{align}
    &W_{0_1}(a,x)=W_{Q_{0_1}}(x_i)|_{x_1=xa^2,\ x_2=x } \nonumber \\ &\text{Li}_2(a^2x)-\text{Li}_2(x)=\text{Li}_2(x_1)-\text{Li}_2(x_2)|_{x_1=xa^2,\ x_2=x }
\end{align}
An interested reader can try to work out the same exercise for the unreduced trefoil \cite{ekholm2020physicsgeometryknotsquiverscorrespondence}. 
\bibliographystyle{utphys}
\bibliography{ref}

\providecommand{\href}[2]{#2}\begingroup\raggedright\begin{thebibliography}{10}

\bibitem{Witten:1988hf}
E.~Witten, ``{Quantum Field Theory and the Jones Polynomial},'' \href{http://dx.doi.org/10.1007/BF01217730}{{\em Commun. Math. Phys.} {\bfseries 121} (1989) 351--399}.

\bibitem{khovanov1999categorificationjonespolynomial}
M.~Khovanov, ``A categorification of the jones polynomial,'' 1999.
\newblock \url{https://arxiv.org/abs/math/9908171}.

\bibitem{dunfield2005superpolynomialknothomologies}
N.~M. Dunfield, S.~Gukov, and J.~Rasmussen, ``The superpolynomial for knot homologies,'' 2005.
\newblock \url{https://arxiv.org/abs/math/0505662}.

\bibitem{Kucharski_2017}
P.~Kucharski, M.~Reineke, M.~Stošić, and P.~Sułkowski, ``Bps states, knots, and quivers,'' \href{http://dx.doi.org/10.1103/physrevd.96.121902}{{\em Physical Review D} {\bfseries 96} no.~12, (Dec., 2017) }. \url{http://dx.doi.org/10.1103/PhysRevD.96.121902}.

\bibitem{Kucharski_2019}
P.~Kucharski, M.~Reineke, M.~Stošić, and P.~Sułkowski, ``Knots-quivers correspondence,'' \href{http://dx.doi.org/10.4310/atmp.2019.v23.n7.a4}{{\em Advances in Theoretical and Mathematical Physics} {\bfseries 23} no.~7, (2019) 1849–1902}. \url{http://dx.doi.org/10.4310/ATMP.2019.v23.n7.a4}.

\bibitem{gopakumar1998mtheorytopologicalstringsii}
R.~Gopakumar and C.~Vafa, ``M-theory and topological strings--ii,'' 1998.
\newblock \url{https://arxiv.org/abs/hep-th/9812127}.

\bibitem{Gopakumar:1998ii}
R.~Gopakumar and C.~Vafa, ``{M theory and topological strings. 1.},'' \href{http://arxiv.org/abs/hep-th/9809187}{{\ttfamily arXiv:hep-th/9809187}}.

\bibitem{Ooguri_2000}
H.~Ooguri and C.~Vafa, ``Knot invariants and topological strings,'' \href{http://dx.doi.org/10.1016/s0550-3213(00)00118-8}{{\em Nuclear Physics B} {\bfseries 577} no.~3, (June, 2000) 419–438}. \url{http://dx.doi.org/10.1016/S0550-3213(00)00118-8}.

\bibitem{Labastida_2001}
J.~M.~F. Labastida and M.~Mariño, ``Polynomial invariants for torus knots¶and topological strings,'' \href{http://dx.doi.org/10.1007/s002200100374}{{\em Communications in Mathematical Physics} {\bfseries 217} no.~2, (Mar., 2001) 423–449}. \url{http://dx.doi.org/10.1007/s002200100374}.

\bibitem{Labastida_2000}
J.~M. Labastida, M.~Mariño, and C.~Vafa, ``Knots, links and branes at large n,'' \href{http://dx.doi.org/10.1088/1126-6708/2000/11/007}{{\em Journal of High Energy Physics} {\bfseries 2000} no.~11, (Nov., 2000) 007–007}. \url{http://dx.doi.org/10.1088/1126-6708/2000/11/007}.

\bibitem{ekholm2020physicsgeometryknotsquiverscorrespondence}
T.~Ekholm, P.~Kucharski, and P.~Longhi, ``Physics and geometry of knots-quivers correspondence,'' 2020.
\newblock \url{https://arxiv.org/abs/1811.03110}.

\bibitem{doi:10.1142/S0218216596000163}
H.~R. MORTON and P.~R. CROMWELL, ``Distinguishing mutants by knot polynomials,'' \href{http://dx.doi.org/10.1142/S0218216596000163}{{\em Journal of Knot Theory and Its Ramifications} {\bfseries 05} no.~02, (1996) 225--238}. \url{https://doi.org/10.1142/S0218216596000163}.

\bibitem{Crane_1994}
L.~Crane and I.~B. Frenkel, ``Four-dimensional topological quantum field theory, hopf categories, and the canonical bases,'' \href{http://dx.doi.org/10.1063/1.530746}{{\em Journal of Mathematical Physics} {\bfseries 35} no.~10, (Oct., 1994) 5136–5154}. \url{http://dx.doi.org/10.1063/1.530746}.

\bibitem{Gukov:2007ck}
S.~Gukov, ``{Gauge theory and knot homologies},'' \href{http://dx.doi.org/10.1002/prop.200610385}{{\em Fortsch. Phys.} {\bfseries 55} (2007) 473--490}, \href{http://arxiv.org/abs/0706.2369}{{\ttfamily arXiv:0706.2369 [hep-th]}}.

\bibitem{Witten:1988ze}
E.~Witten, ``{Topological Quantum Field Theory},'' \href{http://dx.doi.org/10.1007/BF01223371}{{\em Commun. Math. Phys.} {\bfseries 117} (1988) 353}.

\bibitem{witten1994monopolesfourmanifolds}
E.~Witten, ``Monopoles and four-manifolds,'' 1994.
\newblock \url{https://arxiv.org/abs/hep-th/9411102}.

\bibitem{khovanov2004matrixfactorizationslinkhomology}
M.~Khovanov and L.~Rozansky, ``Matrix factorizations and link homology,'' 2004.
\newblock \url{https://arxiv.org/abs/math/0401268}.

\bibitem{Gukov_2005}
S.~Gukov, A.~Schwarz, and C.~Vafa, ``Khovanov-rozansky homology and topological strings,'' \href{http://dx.doi.org/10.1007/s11005-005-0008-8}{{\em Letters in Mathematical Physics} {\bfseries 74} no.~1, (Oct., 2005) 53–74}. \url{http://dx.doi.org/10.1007/s11005-005-0008-8}.

\bibitem{Gukov_2016}
S.~Gukov, S.~Nawata, I.~Saberi, M.~Stošić, and P.~Sułkowski, ``Sequencing bps spectra,'' \href{http://dx.doi.org/10.1007/jhep03(2016)004}{{\em Journal of High Energy Physics} {\bfseries 2016} no.~3, (Mar., 2016) }. \url{http://dx.doi.org/10.1007/JHEP03(2016)004}.

\bibitem{Gukov_2013}
S.~Gukov and M.~Stošić, \href{http://dx.doi.org/10.2140/gtm.2012.18.309}{``Homological algebra of knots and bps states,''} in {\em Proceedings of the Freedman Fest}, vol.~18, p.~309–367.
\newblock Mathematical Sciences Publishers, May, 2013.
\newblock \url{http://dx.doi.org/10.2140/gtm.2012.18.309}.

\bibitem{Taubes_2006}
C.~H. Taubes, \href{http://dx.doi.org/10.2140/gtm.2006.8.73}{``Lagrangians for the gopakumar–vafa conjecture,''} in {\em The interaction of finite-type and Gromov--Witten invariants (BIRS 2003)}, p.~73–95.
\newblock Mathematical Sciences Publishers, Apr., 2006.
\newblock \url{http://dx.doi.org/10.2140/gtm.2006.8.73}.

\bibitem{Aganagic:2015tka}
M.~Aganagic, ``{String theory and math: Why this marriage may last. Mathematics and dualities of quantum physics},'' \href{http://dx.doi.org/10.1090/bull/1517}{{\em Bull. Am. Math. Soc.} {\bfseries 53} no.~1, (2016) 93--115}, \href{http://arxiv.org/abs/1508.06642}{{\ttfamily arXiv:1508.06642 [hep-th]}}.

\bibitem{labastida2001newpointviewtheory}
J.~M.~F. Labastida and M.~Marino, ``A new point of view in the theory of knot and link invariants,'' 2001.
\newblock \url{https://arxiv.org/abs/math/0104180}.

\bibitem{Garoufalidis_2016}
S.~Garoufalidis, P.~Kucharski, and P.~Sułkowski, ``Knots, bps states, and algebraic curves,'' \href{http://dx.doi.org/10.1007/s00220-016-2682-z}{{\em Communications in Mathematical Physics} {\bfseries 346} no.~1, (July, 2016) 75–113}. \url{http://dx.doi.org/10.1007/s00220-016-2682-z}.

\bibitem{Cooper1994}
D.~Cooper, M.~Culler, H.~Gillet, D.~D. Long, and P.~B. Shalen, ``Plane curves associated to character varieties of 3-manifolds,'' \href{http://dx.doi.org/10.1007/BF01231526}{{\em Inventiones mathematicae} {\bfseries 118} no.~1, (Dec, 1994) 47--84}. \url{https://doi.org/10.1007/BF01231526}.

\bibitem{aganagic2012largendualitymirror}
M.~Aganagic and C.~Vafa, ``Large n duality, mirror symmetry, and a q-deformed a-polynomial for knots,'' 2012.
\newblock \url{https://arxiv.org/abs/1204.4709}.

\bibitem{Fuji:2013rra}
H.~Fuji and P.~Sulkowski, ``{Super-A-polynomial},'' \href{http://dx.doi.org/10.1090/pspum/090}{{\em Proc. Symp. Pure Math.} {\bfseries 90} (2015) 277--304}, \href{http://arxiv.org/abs/1303.3709}{{\ttfamily arXiv:1303.3709 [math.AG]}}.

\bibitem{Fuji_2013}
H.~Fuji, S.~Gukov, M.~Stošić, and P.~Sulkowski, ``3d analogs of argyres-douglas theories and knot homologies,'' \href{http://dx.doi.org/10.1007/jhep01(2013)175}{{\em Journal of High Energy Physics} {\bfseries 2013} no.~1, (Jan., 2013) }. \url{http://dx.doi.org/10.1007/JHEP01(2013)175}.

\bibitem{Efimov_2012}
A.~I. Efimov, ``Cohomological hall algebra of a symmetric quiver,'' \href{http://dx.doi.org/10.1112/s0010437x12000152}{{\em Compositio Mathematica} {\bfseries 148} no.~4, (May, 2012) 1133–1146}. \url{http://dx.doi.org/10.1112/S0010437X12000152}.

\bibitem{kontsevich2008stabilitystructuresmotivicdonaldsonthomas}
M.~Kontsevich and Y.~Soibelman, ``Stability structures, motivic donaldson-thomas invariants and cluster transformations,'' 2008.
\newblock \url{https://arxiv.org/abs/0811.2435}.

\bibitem{10.1093/imrn/rnad033}
M.~Reineke, B.~Rhoades, and V.~Tewari, ``Zonotopal algebras, orbit harmonics, and donaldson–thomas invariants of symmetric quivers,'' \href{http://dx.doi.org/10.1093/imrn/rnad033}{{\em International Mathematics Research Notices} {\bfseries 2023} no.~23, (03, 2023) 20169--20210}, \href{http://arxiv.org/abs/https://academic.oup.com/imrn/article-pdf/2023/23/20169/54019853/rnad033.pdf}{{\ttfamily https://academic.oup.com/imrn/article-pdf/2023/23/20169/54019853/rnad033.pdf}}. \url{https://doi.org/10.1093/imrn/rnad033}.

\bibitem{Jankowski_2023}
J.~Jankowski, P.~Kucharski, H.~Larraguível, D.~Noshchenko, and P.~Sułkowski, ``Quiver diagonalization and open bps states,'' \href{http://dx.doi.org/10.1007/s00220-023-04753-2}{{\em Communications in Mathematical Physics} {\bfseries 402} no.~2, (June, 2023) 1551–1584}. \url{http://dx.doi.org/10.1007/s00220-023-04753-2}.

\bibitem{Ekholm:2018eee}
T.~Ekholm, P.~Kucharski, and P.~Longhi, ``{Physics and geometry of knots-quivers correspondence},'' \href{http://dx.doi.org/10.1007/s00220-020-03840-y}{{\em Commun. Math. Phys.} {\bfseries 379} no.~2, (2020) 361--415}, \href{http://arxiv.org/abs/1811.03110}{{\ttfamily arXiv:1811.03110 [hep-th]}}.

\bibitem{meinhardt2016donaldsonthomasinvariantsversusintersection}
S.~Meinhardt and M.~Reineke, ``Donaldson-thomas invariants versus intersection cohomology of quiver moduli,'' 2016.
\newblock \url{https://arxiv.org/abs/1411.4062}.

\bibitem{Franzen_2018}
H.~Franzen and M.~Reineke, ``Semistable chow–hall algebras of quivers and quantized donaldson–thomas invariants,'' \href{http://dx.doi.org/10.2140/ant.2018.12.1001}{{\em Algebra \& Number Theory} {\bfseries 12} no.~5, (July, 2018) 1001–1025}. \url{http://dx.doi.org/10.2140/ant.2018.12.1001}.

\bibitem{gorsky2013quadruplygradedcoloredhomologyknots}
E.~Gorsky, S.~Gukov, and M.~Stosic, ``Quadruply-graded colored homology of knots,'' 2013.
\newblock \url{https://arxiv.org/abs/1304.3481}.

\bibitem{Kucharski_2016}
P.~Kucharski and P.~Sułkowski, ``Bps counting for knots and combinatorics on words,'' \href{http://dx.doi.org/10.1007/jhep11(2016)120}{{\em Journal of High Energy Physics} {\bfseries 2016} no.~11, (Nov., 2016) }. \url{http://dx.doi.org/10.1007/JHEP11(2016)120}.

\bibitem{luo2016integralitystructurestopologicalstrings}
W.~Luo and S.~Zhu, ``Integrality structures in topological strings i: framed unknot,'' 2016.
\newblock \url{https://arxiv.org/abs/1611.06506}.

\bibitem{zhu2018topologicalstringsquivervarieties}
S.~Zhu, ``Topological strings, quiver varieties and rogers-ramanujan identities,'' 2018.
\newblock \url{https://arxiv.org/abs/1707.00831}.

\bibitem{Kucharski:2023jds}
P.~Kucharski, H.~Larragu\i{}vel, D.~Noshchenko, and P.~Su\l{}kowski, ``{Unlinking symmetric quivers},'' \href{http://arxiv.org/abs/2312.14905}{{\ttfamily arXiv:2312.14905 [hep-th]}}.

\bibitem{Jankowski:2021flt}
J.~Jankowski, P.~Kucharski, H.~Larragu\'\i{}vel, D.~Noshchenko, and P.~Su\l{}kowski, ``{Permutohedra for knots and quivers},'' \href{http://dx.doi.org/10.1103/PhysRevD.104.086017}{{\em Phys. Rev. D} {\bfseries 104} no.~8, (2021) 086017}, \href{http://arxiv.org/abs/2105.11806}{{\ttfamily arXiv:2105.11806 [hep-th]}}.

\bibitem{Vafa_2014}
C.~Vafa, M.~Aganagic, T.~Ekholm, and L.~Ng, ``Topological strings, d-model, and knot contact homology,'' \href{http://dx.doi.org/10.4310/atmp.2014.v18.n4.a3}{{\em Advances in Theoretical and Mathematical Physics} {\bfseries 18} no.~4, (2014) 827–956}. \url{http://dx.doi.org/10.4310/ATMP.2014.v18.n4.a3}.

\bibitem{ekholm2025skeinsbranes}
T.~Ekholm and V.~Shende, ``Skeins on branes,'' 2025.
\newblock \url{https://arxiv.org/abs/1901.08027}.

\bibitem{Dimofte_2011}
T.~Dimofte, S.~Gukov, and L.~Hollands, ``Vortex counting and lagrangian 3-manifolds,'' \href{http://dx.doi.org/10.1007/s11005-011-0531-8}{{\em Letters in Mathematical Physics} {\bfseries 98} no.~3, (Oct., 2011) 225–287}. \url{http://dx.doi.org/10.1007/s11005-011-0531-8}.

\bibitem{Terashima_2011}
Y.~Terashima and M.~Yamazaki, ``$ {\text{sl} }\left( {2,\mathbb{R}} \right) $ chern-simons, liouville, and gauge theory on duality walls,'' \href{http://dx.doi.org/10.1007/jhep08(2011)135}{{\em Journal of High Energy Physics} {\bfseries 2011} no.~8, (Aug., 2011) }. \url{http://dx.doi.org/10.1007/JHEP08(2011)135}.

\bibitem{Yagi_2013}
J.~Yagi, ``3d tqft from 6d scft,'' \href{http://dx.doi.org/10.1007/jhep08(2013)017}{{\em Journal of High Energy Physics} {\bfseries 2013} no.~8, (Aug., 2013) }. \url{http://dx.doi.org/10.1007/JHEP08(2013)017}.

\bibitem{Lee_2013}
S.~Lee and M.~Yamazaki, ``3d chern-simons theory from m5-branes,'' \href{http://dx.doi.org/10.1007/jhep12(2013)035}{{\em Journal of High Energy Physics} {\bfseries 2013} no.~12, (Dec., 2013) }. \url{http://dx.doi.org/10.1007/JHEP12(2013)035}.

\bibitem{cordova2013complexchernsimonsm5branessquashed}
C.~Cordova and D.~L. Jafferis, ``Complex chern-simons from m5-branes on the squashed three-sphere,'' 2013.
\newblock \url{https://arxiv.org/abs/1305.2891}.

\bibitem{Hwang_2017}
C.~Hwang, P.~Yi, and Y.~Yoshida, ``Fundamental vortices, wall-crossing, and particle-vortex duality,'' \href{http://dx.doi.org/10.1007/jhep05(2017)099}{{\em Journal of High Energy Physics} {\bfseries 2017} no.~5, (May, 2017) }. \url{http://dx.doi.org/10.1007/JHEP05(2017)099}.

\bibitem{ekholm2020highergenusknotcontact}
T.~Ekholm and L.~Ng, ``Higher genus knot contact homology and recursion for colored homfly-pt polynomials,'' 2020.
\newblock \url{https://arxiv.org/abs/1803.04011}.

\bibitem{aganagic2000mirrorsymmetrydbranescounting}
M.~Aganagic and C.~Vafa, ``Mirror symmetry, d-branes and counting holomorphic discs,'' 2000.
\newblock \url{https://arxiv.org/abs/hep-th/0012041}.

\bibitem{Aganagic_2002}
M.~Aganagic, A.~Klemm, and C.~Vafa, ``Disk instantons, mirror symmetry and the duality web,'' \href{http://dx.doi.org/10.1515/zna-2002-1-201}{{\em Zeitschrift für Naturforschung A} {\bfseries 57} no.~1–2, (Feb., 2002) 1–28}. \url{http://dx.doi.org/10.1515/zna-2002-1-201}.

\bibitem{Panfil_2018}
M.~Panfil, M.~Stošić, and P.~Sułkowski, ``Donaldson-thomas invariants, torus knots, and lattice paths,'' \href{http://dx.doi.org/10.1103/physrevd.98.026022}{{\em Physical Review D} {\bfseries 98} no.~2, (July, 2018) }. \url{http://dx.doi.org/10.1103/PhysRevD.98.026022}.

\bibitem{Panfil_2019}
M.~Panfil and P.~Sułkowski, ``Topological strings, strips and quivers,'' \href{http://dx.doi.org/10.1007/jhep01(2019)124}{{\em Journal of High Energy Physics} {\bfseries 2019} no.~1, (Jan., 2019) }. \url{http://dx.doi.org/10.1007/JHEP01(2019)124}.

\bibitem{ekholm2022knothomologiesgeneralizedquiver}
T.~Ekholm, P.~Kucharski, and P.~Longhi, ``Knot homologies and generalized quiver partition functions,'' 2022.
\newblock \url{https://arxiv.org/abs/2108.12645}.

\bibitem{Beem_2014}
C.~Beem, T.~Dimofte, and S.~Pasquetti, ``Holomorphic blocks in three dimensions,'' \href{http://dx.doi.org/10.1007/jhep12(2014)177}{{\em Journal of High Energy Physics} {\bfseries 2014} no.~12, (Dec., 2014) }. \url{http://dx.doi.org/10.1007/JHEP12(2014)177}.

\bibitem{hwang2013factorization3dsuperconformalindex}
C.~Hwang, H.-C. Kim, and J.~Park, ``Factorization of the 3d superconformal index,'' 2013.
\newblock \url{https://arxiv.org/abs/1211.6023}.

\bibitem{Shadchin_2007}
S.~Shadchin, ``On f-term contribution to effective action,'' \href{http://dx.doi.org/10.1088/1126-6708/2007/08/052}{{\em Journal of High Energy Physics} {\bfseries 2007} no.~08, (Aug., 2007) 052–052}. \url{http://dx.doi.org/10.1088/1126-6708/2007/08/052}.

\end{thebibliography}\endgroup
\end{document}